\documentclass{amsart}

\usepackage{graphicx,epsfig}

\usepackage{amsmath,amssymb,latexsym, amsfonts, amscd, amsthm, xy}

\usepackage{url}

\usepackage{hyperref}

\usepackage{float}
\restylefloat{table}

\input{xy}
\xyoption{all}

\usepackage[usenames]{color}

\makeindex 

=630pt \headheight = 12pt \headsep = 20pt

\usepackage{cite}

\theoremstyle{plain}
\newtheorem{theorem}{Theorem}[section]

\newtheorem{proposition}[theorem]{Proposition}
\newtheorem{corollary}[theorem]{Corollary}

\newtheorem{lemma}[theorem]{Lemma}

\theoremstyle{definition}
\newtheorem{definition}[theorem]{Definition}
\newtheorem{remark}[theorem]{Remark}

\def\bF{\mathbb{F}}

\def\bG{\mathbb{G}}

\def\bZ{\mathbb{Z}}
\def\bR{\mathbb{R}}

\def\cX{\mathcal{X}}

\def\cZ{\mathcal{Z}}

\def\fq{\mathfrak{q}}

\def\deg{\mathrm{deg}}

\begin{document}

\title{Function field analogues of Bang--Zsigmondy's theorem and Feit's theorem}

\author{Nguyen Ngoc Dong Quan}

\date{August 30, 2015}

\address{Department of Mathematics \\
         The University of Texas at Austin \\
         Austin, TX 78712 \\
         USA}

\email{\href{mailto:dongquan.ngoc.nguyen@gmail.com}{\tt dongquan.ngoc.nguyen@gmail.com}}

\maketitle

\tableofcontents

\section{Introduction}
\label{Section-Introduction}

In the number field setting, Bang--Zsigmondy's theorem \cite{Bang} \cite{Zsigmondy} states that for any integers $u, m > 1$, there exists a prime divisor $\wp$ of $u^m - 1$ such that $\wp$ does not divide $u^n - 1$ for every integer $n$ with $0 < n < m$, except exactly in the following cases:
\begin{itemize}

\item [(i)] $m = 2$, and $u = 2^s - 1$ for some integer $s \ge 2$; or

\item [(ii)] $m = 6$, and $u = 2$.

\end{itemize}

A prime $\wp$ satisfying the conditions in Bang--Zsigmondy's theorem is called a \textit{Zsigmondy prime} for $(u, m)$. If $\wp$ is a Zsigmondy prime for $(u, m)$ for some integers $u, m > 1$, then the multiplicative order of $u$ modulo $\wp$ is exactly $m$. Bang--Zsigmondy's theorem has many applications; for example, the existence of Zsigmondy primes was used in the original proof of Wedderburn's theorem \cite{Wedderburn}. See also \cite{Artin} for applications of Zsigmondy primes in theory of finite groups.

Feit \cite{Feit-PAMS-1988} observed that if $\wp$ is a Zsigmondy prime for $(u, m)$, then $\wp \equiv 1 \pmod{m}$ since the multiplicative order of $u$ modulo $\wp$ is exactly $m$. The last congruence implies that $\wp \ge m + 1$, which in turn motivated Feit to introduce the following notion of a large Zsigmondy prime: \textit{a prime $\wp$ is called a large Zsigmondy prime for $(u, m)$ if $\wp$ is a Zsigmondy prime for $(u, m)$ such that either $\wp > m + 1$ or $\wp^2$ divides $u^m - 1$.}

In \cite{Feit-PAMS-1988}, Feit proved a refinement of Bang--Zsigmondy's theorem. He showed that for any integers $u, m > 1$, there exists a large Zsigmondy prime for $(u, m)$ except exactly in the following cases:
\begin{itemize}

\item [(i)] $m = 2$ and $u = 2^s3^t - 1$ for some positive integer $s$, and either $t = 0$ or $1$.

\item [(ii)] $u = 2$ and $m = 4, 6, 10, 12, 18$.

\item [(iii)] $u = 3$ and $m = 4, 6$.

\item [(iv)] $u = 5$ and $m = 6$. 

\end{itemize}

It is obvious that Bang--Zsigmondy's theorem follows immediately from Feit's theorem. 

There are many strong analogies \cite{Goss} \cite{Rosen} \cite{Thakur-book} between number fields and function fields. It is well-known (see  \cite[Chapter {\bf3}]{Goss}) that the Carlitz module is an analogue of the multiplicative group $\bG_m$. The aim of this paper is to search for new analogous phenomena between the Carlitz module and the multiplicative group $\bG_m$. We will introduce notions of Zsigmondy primes and large Zsigmondy primes for the Carlitz module, and prove a Carlitz module analogue of Bang--Zsigmondy's theorem and an analogue of Feit's theorem in the Carlitz module context. 

Throughout the paper, let $q = p^s$, where $p$ is a prime and $s$ is a positive integer. Let $\bF_q$ be the finite field of $q$ elements. Let $A = \bF_q[T]$, and let $k = \bF_q(T)$. Let $\bar{k}$ denote an algebraic closure of $k$. Let $\tau$ be the mapping defined by $\tau(x) = x^q$, and let $k\langle \tau \rangle$ denote the twisted polynomial ring. Let $C : A \rightarrow k\langle \tau \rangle$ ($a \mapsto C_a$) be the Carlitz module, namely, $C$ is an $\bF_q$-algebra homomorphism such that $C_T = T + \tau$.

\subsection{A Carlitz module analogue of Bang--Zsigmondy's theorem}

The main ingredient in the notion of Zsigmondy primes in the number field context is the notion of the multiplicative order of an integer modulo a prime. Hence in order to define a Carlitz module analogue of Zsigmondy primes, we need to find a function field replacement for the notion of the multiplicative order of an integer modulo a prime.

In \cite{dongquan-JNT-145-2014}, we introduced the notion of the Carlitz annihilator of a monic prime $\wp \in A$, which is a function field replacement for the multiplicative order of $2$ modulo a prime. In Section \ref{Section-A-Carlitz-module-analogue-of-Zsigmondy-primes} of the present paper, we will introduce a generalization of the Carlitz annihilator of a monic prime to a couple $(u, \wp)$, where $u$ is a nonzero polynomial in $A$ and $\wp$ is a monic prime. If $\wp$ does not divide $u$, the Carlitz annihilator of a couple $(u, \wp)$, denoted by $P_{u, \wp}$, is the unique monic polynomial in $A$ of least positive degree such that $C_{P_{u, \wp}}(u) \equiv 0 \pmod{\wp}$, i.e., $P_{u, \wp}$ divides $m$ for any nonzero polynomial $m \in A$ with $C_m(u) \equiv 0 \pmod{\wp}$. If $\wp$ divides $u$, we simply let $P_{u, \wp} = 1$. 

The Carlitz annihilator of a couple $(u, \wp)$ can be viewed as a replacement in the function field setting for the multiplicative order of an integer $u$ modulo a prime $\wp$ in the number field setting. This analogy is crucial throughout the present work. The basic analogy between the multiplicative group $\bG_m$ and the Carlitz module $C$ that will be used throughout this paper is illustrated in Table \ref{Table-The-analogy-between-the-multiplicative-group-G_m-and-the-Carlitz-module-C}.
\begin{table}[H]
\caption{The analogy between the multiplicative group $\bG_m$ and the Carlitz module $C$}
\label{Table-The-analogy-between-the-multiplicative-group-G_m-and-the-Carlitz-module-C}
\begin{center}
\begin{tabular}{|l|c|r|}

\hline

\textbf{The multiplicative group $\bG_m$} & \textbf{The Carlitz module $C$} \\

$$ & $$ \\

\hline

$u^m - 1$ for $u, m \in \bZ_{> 0}$ & $C_m(u)$ for $m, u \in A$ \\

$$ & $$ \\

\hline

The multiplicative order of $u$ modulo $\wp$  & The Carlitz annihilator of $(u, \wp)$  \\
$$ & $$\\

\hline

\end{tabular}
\end{center}
\end{table}

In the number field context, recall that a prime $\wp$ is a Zsigmondy prime for $(u, m)$ if the multiplicative order of $u$ modulo $\wp$ is exactly $m$. The analogy in Table \ref{Table-The-analogy-between-the-multiplicative-group-G_m-and-the-Carlitz-module-C} suggests that one can define a Carlitz module analogue of Zsigmondy primes as follows. \textit{For any nonzero polynomials $u, m \in A$, a monic prime $\wp$ is called a Zsigmondy prime for $(u, m)$ if the Carlitz annihilator of $(u, \wp)$ is $m$.} In more concrete terms, this means that $\wp$ divides $C_m(u)$, and $C_n(u) \not\equiv 0 \pmod{\wp}$ for any nonzero polynomial $n \in A$ with $\deg(n) < \deg(m)$.

It is natural to ask whether there exists a Zsigmondy prime for a given couple $(u, m)$, where $u, m$ are nonzero polynomials in $A$. Note that if $\deg(m) = \deg(u) = 0$, then $C_m(u) = mu \in \bF_q^{\times}$, and thus there exists no Zsigmondy primes for $(u, m)$. Therefore without loss of generality, one can modify the last question by adding the assumption that at least one of $m, u$ is of positive degree. As will be explained in Section \ref{Section-A-Carlitz-module-analogue-of-Zsigmondy-primes}, it also suffices to assume further that $m, u$ are monic polynomials.

Determining when there exists no Zsigmondy primes for a given couple $(u, m)$ is a Carlitz module analogue of the classical theorem of Bang--Zsigmondy. The first main result in this paper asserts that there exists a Zsigmondy prime for a given couple $(u, m)$, except exactly in some exceptional cases that can be explicitly determined; more precisely, we obtain the following theorem.

\begin{theorem}[See Theorem \ref{Theorem-The-full-general-version-of-Zsigmondy-Theorem-with-p-not--equal-to-2} and Theorem \ref{Theorem-Zsigmondy-theorem-in-characteristic-two}]
\label{Theorem-The-1st-main-theorem-in-my-paper-summerized-in-the-introduction}

Assume that $q > 2$. Let $m, u$ be monic polynomials in $A$ such that at least one of them is of positive degree. Then there exists a Zsigmondy prime for $(u, m)$ except exactly in the following cases:
\begin{itemize}

\item [(i)] $q = 3$, $u = 1$, and $m = (\wp - 1)\wp$, where $\wp$ is an arbitrary monic prime of degree 1 in $\bF_3[T]$.

\item [(ii)] $q = 2^2$, $u = 1$, and $m = (\wp - 1)\wp$, where $\wp$ is an arbitrary monic prime of degree 1 in $\bF_{2^2}[T]$.

\end{itemize}

\end{theorem}

Theorem \ref{Theorem-The-1st-main-theorem-in-my-paper-summerized-in-the-introduction} will be split into two parts. The first part is Theorem \ref{Theorem-The-full-general-version-of-Zsigmondy-Theorem-with-p-not--equal-to-2} that considers the case when $p \ne 2$. The second part is Theorem \ref{Theorem-Zsigmondy-theorem-in-characteristic-two} that treats Theorem \ref{Theorem-The-1st-main-theorem-in-my-paper-summerized-in-the-introduction} in the case when $p = 2$.

A Carlitz module analogue of the classical Bang--Zsigmondy theorem was first considered by Bae in \cite{Bae}. In fact, Bae considered a more general analogue of the classical Bang--Zsigmondy theorem, which can be viewed as a function field analogue of the classical Zsigmondy theorem. (See \cite{Zsigmondy} or \cite{Birkhoff-Vandiver} for an account of the classical Zsigmondy theorem. Note that Birkhoff and Vandiver \cite{Birkhoff-Vandiver} independently discovered Zsigmondy's theorem in 1904 after the theorem was first proved by Zsigmondy \cite{Zsigmondy} in 1892.) The main result in Bae \cite{Bae} (see \cite[Theorem 4.10]{Bae}) is erroneous. Before pointing out the error in \cite[Theorem 4.10]{Bae}, let us recall the statement of Theorem 4.10 in \cite{Bae}.

Let $u, v$ be relatively prime elements in $A$, and let $m$ be a monic polynomial in $A$. Set 
\begin{align*}
\cZ_m(u, v) = v^{q^{\deg(m)}}C_m(u/v).
\end{align*}
Following Bae \cite{Bae}, a monic prime $\wp$ is called a \textit{primitive factor} of $\cZ_m(u, v)$ if $\cZ_m(u, v) \equiv 0 \pmod{\wp}$, and $\cZ_n(u, v) \not\equiv 0 \pmod{\wp}$ for any monic divisor $n$ of $m$ with $n \ne m$. Note that when $v = 1$, the notion of primitive factors in \cite{Bae} agrees with that of Zsigmondy primes introduced in this paper. 

Theorem 4.10 in Bae \cite{Bae} claims that if $q > 2$ and $\deg(m) > 0$, then $\cZ_m(u, v)$ has at least one primitive factor except exactly in the following case:
\begin{itemize}
 
\item [(PF)] $q = 3$, $u = \pm 1$, $v = \pm 1$, and $m = (\wp - 1)\wp$, where $\wp$ is an arbitrary monic prime of degree one in $\bF_3[T]$.

\end{itemize}
We now provide a counterexample to Theorem 4.10 in Bae \cite{Bae}. Indeed, if one takes $q = 2^2$, $u = v = 1$, and $m = (\wp - 1)\wp$, where $\wp$ is an arbitrary monic prime of degree one in $\bF_{2^2}[T]$, then Table \ref{Table-There-are-no-Zsigmondy-primes-for-(1,m)-when-q-equals-4} in the proof of Lemma \ref{Lemma-The-2nd-lemma-for-Zsigmondy-theorem-in-char-2} shows that there exits no Zsigmondy prime for $(u, m)$, i.e., $\cZ_m(u, v)$ has no primitive factors in this case, which provides a counterexample to Theorem 4.10 in Bae \cite{Bae}.

We should also note that in his Ph.D. thesis (see \cite[Theorem 4.2.10]{Bamunoba}), Bamunoba, following very closely the techniques in Bae \cite{Bae}, attempted to give a different proof of Theorem 4.10 in Bae \cite{Bae}. Due to the counterexample to Theorem 4.10 in Bae \cite{Bae} that we pointed out above, the proof of Theorem 4.2.10 in Bamunoba \cite{Bamunoba} is also erroneous. 

Note that the techniques that Bae exploited in \cite{Bae} are based on the work of Birkhoff and Vandiver \cite{Birkhoff-Vandiver}. We, however, use completely different arguments from Bae \cite{Bae} to prove Theorem \ref{Theorem-The-1st-main-theorem-in-my-paper-summerized-in-the-introduction}, and the strategy of our proof is similar to the work of Roitman \cite{Roitman}. Let us now describe the strategy of our proof of Theorem \ref{Theorem-The-1st-main-theorem-in-my-paper-summerized-in-the-introduction} in detail.

In Section \ref{Section-A-Carlitz-module-analogue-of-Zsigmondy-primes}, we obtain a function field analogue of L\"uneburg's theorem (see L\"uneburg \cite[Satz 1]{Luneburg} or Roitman \cite[Proposition 2]{Roitman} for an account of L\"uneburg's theorem in the number field context) that describes a sufficient and necessary condition under which a prime $\wp$ is a non-Zsigmondy prime for $(u, m)$. 

\begin{theorem}
\label{Thm-Luneburg-theorem-in-the-introduction}
$(\text{See Theorem \ref{Theorem-The-1st-version-of-Zsigmondy-Theorem} or Corollary \ref{C-Luneburg-theorem-when-q->-2}})$
  
Assume that $q > 2$. Let $m, u$ be monic polynomials in $A$ such that $m$ is of positive degree. Let $\wp$ be a monic prime dividing $\Psi_m(u)$, where $\Psi_m(x) \in A[x]$ denotes the $m$-th cyclotomic polynomial (that will be reviewed in Section \ref{Section-Cyclotomic-polynomials-over-function-fields}). Let $P_{u, \wp}$ be the Carlitz annihilator of $(u, \wp)$. Then
\begin{itemize}

\item [(i)] $\wp$ is a non-Zsigmondy prime for $(u, m)$ if and only if $\wp$ divides $m$.

\item [(ii)] if $\wp$ is a non-Zsigmondy prime for $(u, m)$, then $m = P_{u, \wp} \wp^s$ for some positive integer $s$. Furthermore $\wp^2$ does not divide $\Psi_m(u)$.

\end{itemize}

\end{theorem}

\begin{remark}

Note that Bae (see \cite[Proposition 4.4]{Bae}) also obtained a more general version of part $(i)$ of Theorem \ref{Thm-Luneburg-theorem-in-the-introduction}.

In his Ph.D. thesis, Bamunoba (see \cite[Corollary 4.1.5]{Bamunoba}) independently obtained part $(i)$ of Theorem \ref{Thm-Luneburg-theorem-in-the-introduction} with a slightly different proof. Bamunoba (see \cite[Lemma 4.1.6]{Bamunoba}) also independently proved part $(ii)$ of Theorem \ref{Thm-Luneburg-theorem-in-the-introduction} under the more restrictive assumption that $\gcd(u, m) = 1$. 
\end{remark}

Using Theorem \ref{Thm-Luneburg-theorem-in-the-introduction}, and under the assumption that there are no Zsigmondy primes for a pair $(u, m)$, we deduce that the prime factorization of $\Psi_m(u)$ in $A$ is of very special form. If the characteristic $p$ of $k$ is not equal to 2, then $\Psi_m(u)$ is a prime dividing $m$ (see Corollary \ref{Corollary-The-description-of-non-Zsigmondy-primes}). If $p = 2$, then either $\Psi_m(u)$ is a prime dividing $m$ or $\Psi_m(u) = \epsilon \wp(\wp - 1)$, where $\epsilon \in \bF_q^{\times}$ and both $\wp$ and $\wp - 1$ are monic prime divisors of $m$ (see Lemma \ref{Lemma-The-prime-factorization-of-Psi-m-u-when-there-are-no-Zsigmondy-primes-for-u-m-and-p=2}). In either case, by deriving several lower bounds for $\deg(\Psi_m(u))$, we show that except the exceptional cases $(i), (ii)$ in Theorem \ref{Theorem-The-1st-main-theorem-in-my-paper-summerized-in-the-introduction}, the prime factorization of $\Psi_m(u)$ can not fall into these forms, and of course this proves that there must exist a Zsigmondy prime for $(u, m)$.

\subsection{A Carlitz module analogue of Feit's theorem}

Before stating the second main result in this paper, let us describe, for each monic prime $\wp$ in $A$, another $A$-module structure of $A/\wp A$ from which a function field analogue of Fermat's little theorem immediately follows.

Let $\wp$ be a monic prime in $A$. For each $n \in A$, denote by $\bar{n}$ the image of $n$ in $A/\wp A$. We denote by $\bar{C}$ the reduction of $C$ modulo $\wp$. The action of $\bar{C}$ on $A/\wp A$ is given by $\bar{C}_T(u) = \bar{T}u + u^q$ for each $u \in A/\wp A$. Under the action of $\bar{C}$, one obtains another $A$-module structure of $A/\wp A$, which we will denote by $(A/\wp A)_C$. It is known (see Hsu \cite[page 249]{Hsu}) that $(A/\wp A)_C$ is isomorphic to $A/(\wp - 1)A$. From this isomorphism, one immediately obtains the following result which can be viewed as a function field analogue of Fermat's little theorem. (See also Hayes \cite[Proposition 2.4]{Hayes} for another proof.)

\begin{lemma}
\label{Lemma-An-analogue-of-Fermat-little-theorem}

Let $u$ be a nonzero polynomial in $A$, and let $\wp$ be a monic prime in $A$. Then
\begin{align*}
C_{\wp - 1}(u) \equiv 0 \pmod{\wp}.
\end{align*}

\end{lemma}

Now let $\wp$ be a Zsigmondy prime for $(u, m)$, where $u, m$ are monic polynomials in $A$. Then the Carlitz annihilator of $(u, \wp)$ is $m$, i.e., $m$ is the unique monic polynomial of least positive degree such that $C_m(u) \equiv 0 \pmod{\wp}$. Hence it follows from Lemma \ref{Lemma-An-analogue-of-Fermat-little-theorem} that $m$ divides $\wp - 1$, which implies that $\deg(\wp) = \deg(\wp - 1) \ge \deg(m)$. The last inequality and the notion of large Zsigmondy primes introduced in Feit \cite{Feit-PAMS-1988} motivate a notion of a function field analogue of large Zsigmondy primes as follows: \textit{A Zsigmondy prime $\wp$ for $(u, m)$ is called a large Zsigmondy prime for $(u, m)$ if either $\deg(\wp) > \deg(m)$ or $\wp^2$ divides $C_m(u)$.} 

Classifying all pairs $(u, m)$ such that there exists a large Zsigmondy prime for $(u, m)$ can be viewed as a function field analogue of Feit's theorem. We are now in a position to state our most important result in this paper, a function field analogue of Feit's theorem.

\begin{theorem}[See Theorem \ref{Theorem-The-main-theorem-about-large-Zsigmondy-primes}]
\label{Theorem-The-2nd-main-theorem-in-my-paper-summerized-in-the-introduction}

Assume that $q > 2$. Let $m, u$ be monic polynomials in $A$ such that at least one of them is of positive degree. Then there exists a large Zsigmondy prime for $(u, m)$ except in some exceptional cases that can be explicitly determined. (Theorem \ref{Theorem-The-main-theorem-about-large-Zsigmondy-primes} explicitly lists all triples $(q, u, m)$ in the exceptional cases.)

\end{theorem}

We will prove Theorem \ref{Theorem-The-2nd-main-theorem-in-my-paper-summerized-in-the-introduction} in Section \ref{Section-Large-Zsigmondy-primes}. It is obvious that Theorem \ref{Theorem-The-1st-main-theorem-in-my-paper-summerized-in-the-introduction} follows immediately from Theorem \ref{Theorem-The-2nd-main-theorem-in-my-paper-summerized-in-the-introduction}.

\subsection{Notation}
\label{Notation}

Every nonzero element $m \in A$ is of the form $m = \alpha_n T^n + \cdots + \alpha_1 T + \alpha_0$, where the $\alpha_i$ are elements in $\bF_q$ and $\alpha_n \ne 0$. When $m$ is of the form as above, we say that the degree of $m$ is $n$. In notation, we write $\deg(m) = n$. We use the standard convention that $\deg(0) = -\infty$. With this convention, one obtains the degree function $\deg : A \rightarrow \bZ \cup \{-\infty\}$ in an obvious way. 

With $m$ of the above form, we say that the \textit{leading coefficient} of $m$ is $\alpha_n$. 

For a polynomial $m \in A$ of positive degree, we define $\Phi(m)$ to be the number of nonzero polynomials of degree less than $\deg(m)$ and relatively prime to $m$. The function $\Phi(\cdot)$ is a function field analogue of the classical Euler $\phi$-function.

Let $m = \alpha \wp_1^{s_1}\cdots \wp_h^{s_h}$ be the prime factorization of $m$, where $\alpha \in \bF_q^{\times}$, the $\wp_i$ are monic primes in $A$, and the $s_i$ are positive integers. It is well-known (see \cite[Proposition 1.7]{Rosen}) that
\begin{align*}
\Phi(m) = \prod_{i = 1}^h\Phi(\wp_i^{s_i}) = \prod_{i = 1}^h(q^{\deg(\wp_i^{s_i})} - q^{\deg(\wp_i^{s_i - 1})}).
\end{align*}
In particular, this implies that when $m = \wp^s$ for some monic prime $\wp$ and some positive integer $s$, then
\begin{align*}
\Phi(\wp^s) = q^{\deg(\wp^{s})} - q^{\deg(\wp^{s - 1})}.
\end{align*}

\section{Cyclotomic polynomials over function fields}
\label{Section-Cyclotomic-polynomials-over-function-fields}

In this section, we prove some results about cyclotomic polynomials over function fields that will be useful in subsequent sections. We begin by recalling the definition of cyclotomic polynomials in function fields, and some well-known results of cyclotomic polynomials in function fields whose proofs can be found, for example, in \cite{Bae-Hahn-BKMS-1992} or \cite{Villa-Salvador}.

Let $m$ be a polynomial of positive degree, and set $\Lambda_m := \{\lambda \in \bar{k}\; | \; C_m(\lambda) = 0 \}$. We define a \textit{primitive $m$-th root of $C$} to be a root of the polynomial $C_m(x) \in A[x]$ that generates the $A$-module $\Lambda_m$. We fix a primitive $m$-th root of $C$, and denote it by $\lambda_m$. 

Recall that the \textit{$m$-th cyclotomic polynomial}, denoted by $\Psi_m(x)$, is the minimal polynomial of $\lambda_m$ over $k$, i.e., $\Psi_m(x) \in k[x]$ is the monic irreducible polynomial of least degree such that $\Psi_m(\lambda_m) = 0$. It is well-known that $\Psi_m(x) \in A[x]$. When $m = \wp^s$ for some monic prime $\wp$ and some positive integer $s$, we know from \cite[Proposition 2.4]{Hayes} that 
\begin{align}
\label{Eqn-the-cyclotomic-polynomial-at-a-power-of-prime}
\Psi_{\wp^s}(x) = C_{\wp^s}(x)/C_{\wp^{s - 1}}(x).
\end{align}  

The next two results are well-known.

\begin{proposition}
\label{Proposition-The-2nd-basic-result-of-cyclotomic-polynomials}
$(\text{See part (2) in \cite[Proposition 12.3.13]{Villa-Salvador}})$

Let $m$ be a monic polynomial in $A$. Then
\begin{align*}
C_m(x) = \prod_{\substack{b | m, \\ b \; \text{monic}}}\Psi_b(x).
\end{align*}

\end{proposition}

\begin{proposition}
\label{Lemma-The-1st-lemma-about-cyclotomic-polynomials}
$(\text{See \cite[Proposition 1.2(c)]{Bae-Hahn-BKMS-1992}})$

Let $\wp$ be a monic prime in $A$, and let $m$ be a monic polynomial in $A$ such that $\gcd(m, \wp) = 1$. Let $h$ be a positive integer. Then
\begin{align*}
\Psi_{m}(C_{\wp^h}(x)) = \Psi_{m\wp^h}(x)\Psi_m(C_{\wp^{h - 1}}(x)).
\end{align*}

\end{proposition}

Using Proposition \ref{Lemma-The-1st-lemma-about-cyclotomic-polynomials}, we prove the following result that will be useful in the proof of a Carlitz module analogue of L\"uneburg's theorem.

\begin{lemma}
\label{Corollary-The-4th-basic-result-of-cyclotomic-polynomials}

Let $m$ be a monic polynomial in $A$ of positive degree, and let $\fq$ be a monic prime in $A$ such that $\fq$ divides $m$. Then $\Psi_m(x)$ divides $C_m(x)/C_{m/\fq}(x)$ in the polynomial ring $A[x]$.

\end{lemma}

\begin{proof}

We first prove Lemma \ref{Corollary-The-4th-basic-result-of-cyclotomic-polynomials} in the case when $m = \wp^s$ for some monic prime $\wp$ and some positive integer $s$. In this case, we see that $\fq = \wp$, and thus
\begin{align*}
\Psi_m(x) =  \Psi_{\wp^s}(x) = \dfrac{C_{\wp^s}(x)}{C_{\wp^{s-1}}(x)} = \dfrac{C_m(x)}{C_{m/\wp}(x)} = \dfrac{C_m(x)}{C_{m/\fq}(x)},
\end{align*}
which prove Lemma \ref{Corollary-The-4th-basic-result-of-cyclotomic-polynomials} for $m = \wp^s$.

Set $d = \deg(m) \ge 1$. We prove Lemma \ref{Corollary-The-4th-basic-result-of-cyclotomic-polynomials} by induction on $d$.

If $d = 1$, then $m = \wp$ for some monic prime $\wp$ of degree $1$. Thus Lemma \ref{Corollary-The-4th-basic-result-of-cyclotomic-polynomials} holds for $d = 1$.

Assume that Lemma \ref{Corollary-The-4th-basic-result-of-cyclotomic-polynomials} holds for any monic polynomial $m$ of degree less than $d$. We prove that Lemma \ref{Corollary-The-4th-basic-result-of-cyclotomic-polynomials} is true for $d$. Indeed, take any monic polynomial $m$ of degree $d$. If $m$ has exactly one monic prime factor $\wp$, then $m = \wp^s$ for some positive integer $s$, and we already prove that Lemma \ref{Corollary-The-4th-basic-result-of-cyclotomic-polynomials} is true in this case. 

If $m$ has at least two distinct prime factors, then there exists another monic prime $\wp$ with $\wp \ne \fq$. Write $m = n\wp^s$ for some positive integer $s$, where $n$ is a monic polynomial such that $\gcd(n, \wp) = 1$. We know that $\fq$ divides $n$ and $1 \le \deg(n) < \deg(m) = d$, and it thus follows from the induction hypothesis that $\Psi_n(x)$ divides $C_n(x)/C_{n/\fq}(x)$. Hence there exists a polynomial $\Gamma(x) \in A[x]$ such that $\Psi_n(x)\Gamma(x) = \dfrac{C_n(x)}{C_{n/\fq}(x)}.$ Substituting $C_{\wp^s}(x)$ for $x$ in the last equation, we deduce that
\begin{align}
\label{Equation-The-1st-equation-in-the-corollary-about-Psi-m-dividing-C-m(x)-over-C-m/fq(x)}
\Psi_n(C_{\wp^s}(x))\Gamma(C_{\wp^s}(x)) = \dfrac{C_n(C_{\wp^s}(x))}{C_{n/\fq}(C_{\wp^s}(x))} = \dfrac{C_{n\wp^s}(x)}{C_{(n\wp^s)/\fq}(x)} = \dfrac{C_m(x)}{C_{m/\fq}(x)}.
\end{align}

Using equation (\ref{Equation-The-1st-equation-in-the-corollary-about-Psi-m-dividing-C-m(x)-over-C-m/fq(x)}) and applying Proposition \ref{Lemma-The-1st-lemma-about-cyclotomic-polynomials} with $n, \wp^s$ in the roles of $m, \wp^h$, respectively, we deduce that
\begin{align*}
\Psi_m(x)\Psi_n(C_{\wp^{s - 1}}(x))\Gamma(C_{\wp^s}(x)) = \Psi_{n\wp^s}(x)\Psi_n(C_{\wp^{s - 1}}(x)) \Gamma(C_{\wp^s}(x)) = \Psi_n(C_{\wp^s}(x)) \Gamma(C_{\wp^s}(x))= \dfrac{C_m(x)}{C_{m/\fq}(x)},
\end{align*}
which proves Lemma \ref{Corollary-The-4th-basic-result-of-cyclotomic-polynomials}.

\end{proof}

The next two results are well-known, and will be very useful in many places of this paper. For the proof of these results, see Bae \cite[Lemmas 4.6 and 4.8]{Bae}. 
 
\begin{proposition}
\label{Proposition-The-1st-result-about-the-degrees-of-cyclotomic-polynomials}

Let $m$ be an element in $A$, and let $u$ be a polynomial in $A$ of positive degree. Assume that the following condition is true:
\begin{itemize}

\item [(D)] if $\deg(u) = 1$, then $q > 2$.

\end{itemize}
Then 
\begin{align*}
\deg(C_m(u)) =
\begin{cases}
-\infty \; \; &\text{if $m = 0$,} \\
\deg(u)q^{\deg(m)} \; \; &\text{if $m \ne 0$.}
\end{cases}
\end{align*}

\end{proposition}

\begin{proposition}
\label{Proposition-The-2nd-result-about-the-degrees-of-cyclotomic-polynomials-when-u-is-a-unit}

Let $m$ be an element in $A$, and let $u$ be a unit in $\bF_q^{\times}$. Assume that $q > 2$. Then the degree of $C_m(u)$ satisfies
\begin{align*}
\deg(C_m(u)) =
\begin{cases}
-\infty \; &\text{if $m = 0$,} \\
0 \; &\text{if $\deg(m) = 0$}, \\
q^{\deg(m) - 1} \; &\text{otherwise}.
\end{cases}
\end{align*}

\end{proposition}

The following elementary result will be used at many times in proofs of some subsequent results. 

\begin{corollary}
\label{Corollary-C_m(u)-is-nonzero}

Let $m$ be a nonzero polynomial in $A$, and let $u$ be a polynomial in $A$. Assume that the following is true:
\begin{itemize}

\item [$(D*)$] if $\deg(u) = 0$ or $\deg(u) = 1$, then $q > 2$.

\end{itemize}
Then $C_m(u) \ne 0$ if and only if $u \ne 0$.

\end{corollary}

\begin{proof}

If $u = 0$, then $C_m(u) = 0$. If $u \ne 0$, then either $\deg(u) = 0$ or $\deg(u) \ge 1$. If $\deg(u) = 0$, it follows from $(D*)$ that $q > 2$. Since $u \ne 0$, we know that $u$ is a unit in $\bF_q^{\times}$, and it thus follows from Proposition \ref{Proposition-The-2nd-result-about-the-degrees-of-cyclotomic-polynomials-when-u-is-a-unit} that $C_m(u) \ne 0$.

If $\deg(u) \ge 1$, we deduce from $(D*)$ that condition (D) in Proposition \ref{Proposition-The-1st-result-about-the-degrees-of-cyclotomic-polynomials} is satisfied, and it thus follows from Proposition \ref{Proposition-The-1st-result-about-the-degrees-of-cyclotomic-polynomials} that $C_m(u) \ne 0$. 

\end{proof}

In the number field context, Roitman \cite{Roitman} obtained several lower bounds for the value of $\Psi_m(u)$, where $\Psi_m(x) \in \bZ[x]$ is the classical $m$-th cyclotomic polynomial and $u$ is a positive integer. These lower bounds play a significant role in the proofs of the classical Bang--Zsigmondy theorem and the classical Feit theorem that are given in Roitman \cite{Roitman}. In the function field context, in order to measure how \textit{large} a polynomial in $A$ is, one can use the degree function $\deg: A \rightarrow \bZ \cup \{-\infty\}$ in place of the usual absolute value $| \cdot |$ of $\bR$. The aim of the next two lemmas is to compute the value of $\deg(\Psi_m(u))$, where $\Psi_m(x) \in A[x]$ is the $m$-th cyclotomic polynomial over function fields, and $m, u$ are polynomials in $A$. In contrast to the classical case, we can obtain an exact formula for $\deg(\Psi_m(u))$. The next two lemmas will be crucial in the proofs of our main results. 

We first prove an exact formula for $\deg(\Psi_m(u))$ in the case when $u$ is of positive degree.

\begin{lemma}
\label{Lemma-The-degrees-of-cyclotomic-polynomials}

Let $m$ be a monic polynomial in $A$ of positive degree, and let $u$ be a polynomial of positive degree. Assume that (D) in Proposition \ref{Proposition-The-1st-result-about-the-degrees-of-cyclotomic-polynomials} is true. Then
\begin{align*}
\deg(\Psi_m(u)) = \deg(u)\Phi(m),
\end{align*}
where $\Phi(\cdot)$ is the function field analogue of the classical Euler $\phi$-function (see Subsection \ref{Notation} for its definition).

\end{lemma}

\begin{proof}

Let us first consider the case when $m = P^s$ for some monic prime $P \in A$ and $s \in \bZ_{>0}$. By Corollary \ref{Corollary-C_m(u)-is-nonzero}, $C_{P^{s - 1}}(u) \ne 0$, and thus (\ref{Eqn-the-cyclotomic-polynomial-at-a-power-of-prime}) implies that $\Psi_m(u) = \Psi_{P^s}(u) = \dfrac{C_{P^s}(u)}{C_{P^{s-1}}(u)}$. Applying Proposition \ref{Proposition-The-1st-result-about-the-degrees-of-cyclotomic-polynomials}, we see that $\deg(C_{P^{s - 1}}(u)) = \deg(u)q^{\deg(P^{s - 1})}$ and $\deg(C_{P^s}(u)) = \deg(u)q^{\deg(P^s)}$. Thus
\begin{align*}
\deg(\Psi_m(u)) &= \deg(C_{P^s}(u)) - \deg(C_{P^{s - 1}}(u)) = \deg(u)(q^{\deg(P^s)} - q^{\deg(P^{s - 1})}) \\
&= \deg(u)\Phi(P^s) = \deg(u)\Phi(m),
\end{align*}
which proves Lemma \ref{Lemma-The-degrees-of-cyclotomic-polynomials} for $m = P^s$.

Now let $m$ be a monic polynomial in $A$ of positive degree. Write $m$ in the form $m = n\wp^s$ with $s \in \bZ_{>0}$, where $n$ is a monic polynomial and $\wp$ is a monic prime such that $\gcd(n, \wp) = 1$. If $n = 1$, then $m = \wp^s$, and in this case Lemma \ref{Lemma-The-degrees-of-cyclotomic-polynomials} is true as shown above.

If $n$ is of positive degree, then the induction hypothesis tells us that Lemma \ref{Lemma-The-degrees-of-cyclotomic-polynomials} is true for $n$, that is, $\deg(\Psi_n(v)) = \deg(v)\Phi(n)$ for any $v \in A$ of positive degree. Using this fact and Proposition \ref{Proposition-The-1st-result-about-the-degrees-of-cyclotomic-polynomials}, and noting that $\deg(C_{\wp^s}(u)) \ge 1$ and $\deg(C_{\wp^{s - 1}}(u)) \ge 1$, we deduce that
\begin{align*}
\deg(\Psi_n(C_{\wp^s}(u) ) = \deg(C_{\wp^s}(u))\Phi(n) = \deg(u)q^{\deg(\wp^s)}\Phi(n),
\end{align*}
and
\begin{align*}
\deg(\Psi_n(C_{\wp^{s - 1}}(u)) = \deg(C_{\wp^{s - 1}}(u))\Phi(n) = \deg(u)q^{\deg(\wp^{s - 1})}\Phi(n).
\end{align*}
Therefore it follows from Proposition \ref{Lemma-The-1st-lemma-about-cyclotomic-polynomials} that
\begin{align*}
\deg(\Psi_m(u)) &= \deg(\Psi_{n\wp^s}(u)) = \deg(\Psi_n(C_{\wp^s}(u)) - \deg(\Psi_n(C_{\wp^{s - 1}}(u)) \\
&= \deg(u)q^{\deg(\wp^s)}\Phi(n) - \deg(u)q^{\deg(\wp^{s - 1})}\Phi(n) = \deg(u)\Phi(n)(q^{\deg(\wp^s)} - q^{\deg(\wp^{s - 1})}) \\
&= \deg(u)\Phi(n)\Phi(\wp^s) = \deg(u)\Phi(m),
\end{align*}
which proves our contention.

\end{proof}

We now prove a formula for $\deg(\Psi_m(u))$, where $u$ is a unit in $\bF_q^{\times}$.

\begin{lemma}
\label{Lemma-The-degrees-of-cyclotomic-polynomials-when-u-is-a-unit}

Let $m$ be a monic polynomial in $A$ of positive degree, and let $u$ be a unit in $\bF_q^{\times}$. Assume that $q > 2$. Then the degree of $\Psi_m(u)$ satisfies
\begin{align*}
\deg(\Psi_m(u)) =
\begin{cases}
\dfrac{\Phi(m) + (-1)^{h + 1}}{q} \; &\text{if $m = \wp_1\wp_2\cdots\wp_h$, where the $\wp_i$ are distinct monic primes,} \\
\dfrac{\Phi(m)}{q} \; &\text{if there exists a monic prime $\wp$ such that $\wp^2$ divides $m$}.
\end{cases}
\end{align*}

\end{lemma}

\begin{proof}

We consider the cases:

$\star$ \textit{Case 1. $m$ is square-free.} 

One can write $m$ in the form 
\begin{align}
\label{Eqn-of-m-square-fee-in-the-lemma-about-deg-of-CP-with-u-a-unit}
m =  \wp_1\wp_2\cdots\wp_h, 
\end{align}
where the $\wp_i$ are distinct monic primes and $h$ is a positive integer. We prove Lemma \ref{Lemma-The-degrees-of-cyclotomic-polynomials-when-u-is-a-unit} by induction on $h$. If $h = 1$, then $m = \wp$ for some monic prime $\wp$. We know from (\ref{Eqn-the-cyclotomic-polynomial-at-a-power-of-prime}) that $\Psi_m(u) = \dfrac{C_{\wp}(u)}{u}$, and thus Proposition \ref{Proposition-The-2nd-result-about-the-degrees-of-cyclotomic-polynomials-when-u-is-a-unit} implies that $\deg(C_{\wp}(u)) = q^{\deg(\wp) - 1}$. Therefore
\begin{align*}
\deg(\Psi_m(u)) = \deg(C_\wp(u)) - \deg(u) = q^{\deg(\wp) - 1} = \dfrac{\Phi(\wp) + 1}{q} = \dfrac{\Phi(m) + (-1)^{h + 1}}{q},
\end{align*}
which proves Lemma \ref{Lemma-The-degrees-of-cyclotomic-polynomials-when-u-is-a-unit} for $h = 1$.

Assume that Lemma \ref{Lemma-The-degrees-of-cyclotomic-polynomials-when-u-is-a-unit} is true for $h - 1$ with $h \ge 2$. We now prove that Lemma \ref{Lemma-The-degrees-of-cyclotomic-polynomials-when-u-is-a-unit} is true for $h$. Indeed, take any polynomial $m \in A$ of the form (\ref{Eqn-of-m-square-fee-in-the-lemma-about-deg-of-CP-with-u-a-unit}), and set $n = \wp_2\cdots\wp_h$. By Proposition \ref{Proposition-The-2nd-result-about-the-degrees-of-cyclotomic-polynomials-when-u-is-a-unit}, we know that $\deg(C_{\wp_1}(u)) = q^{\deg(\wp_1) - 1} \ge 1$, and it thus follows from Lemma \ref{Lemma-The-degrees-of-cyclotomic-polynomials} that
\begin{align}
\label{Eqn-deg-of-Psi_n(C_wp1(u))-in-lemma-about-deg-of-CP-with-u-a-unit}
\deg(\Psi_n(C_{\wp_1}(u))) = \deg(C_{\wp_1}(u))\Phi(n) = q^{\deg(\wp_1) - 1}\Phi(n).
\end{align}

By the induction hypothesis, we know that
\begin{align}
\label{Eqn-deg-of-Psi_n(u)-in-lemma-about-deg-of-CP-with-u-a-unit}
\deg(\Psi_n(u)) = \deg(\Psi_{\wp_2\cdots\wp_h}(u)) = \dfrac{\Phi(n) + (-1)^h}{q}.
\end{align}
In particular, this implies that $\Psi_n(u) \ne 0$. Applying Proposition \ref{Lemma-The-1st-lemma-about-cyclotomic-polynomials} with $1, \wp_1, n$ in the roles of $h, \wp, m$, respectively, we deduce that
\begin{align*}
\Psi_m(u) = \Psi_{\wp_1n}(u) = \dfrac{\Psi_n(C_{\wp_1}(u))}{\Psi_n(C_{\wp_1^0}(u))} = \dfrac{\Psi_n(C_{\wp_1}(u))}{\Psi_n(C_{1}(u))} = \dfrac{\Psi_n(C_{\wp_1}(u))}{\Psi_n(u)}.
\end{align*}
Thus we deduce from (\ref{Eqn-deg-of-Psi_n(C_wp1(u))-in-lemma-about-deg-of-CP-with-u-a-unit}) and (\ref{Eqn-deg-of-Psi_n(u)-in-lemma-about-deg-of-CP-with-u-a-unit}) that
\begin{align*}
\deg(\Psi_m(u)) &= \deg(\Psi_n(C_{\wp_1}(u))) - \deg(\Psi_n(u)) =  q^{\deg(\wp_1) - 1}\Phi(n) - \dfrac{\Phi(n) + (-1)^h}{q} \\
&= \dfrac{\Phi(n)(q^{\deg(\wp_1)} - 1) + (-1)^{h + 1}}{q} = \dfrac{\Phi(n)\Phi(\wp_1) + (-1)^{h + 1}}{q} = \dfrac{\Phi(m) + (-1)^{h + 1}}{q},
\end{align*}
which proves that Lemma \ref{Lemma-The-degrees-of-cyclotomic-polynomials-when-u-is-a-unit} holds for any monic polynomial $m$ of the form (\ref{Eqn-of-m-square-fee-in-the-lemma-about-deg-of-CP-with-u-a-unit}).

$\star$ \textit{Case 2. $m$ is a monic polynomial such that $\wp^2$ divides $m$ for some monic prime $\wp$.} 

Write $m = n\wp^s$, where $s \in \bZ$ with $s \ge 2$, $\wp$ is a monic prime, and $n$ is a monic polynomial with $\gcd(n, \wp) = 1$. Using the same arguments as in Case 1, one obtains the proof of Case 2.

\end{proof}

\section{A Carlitz module analogue of Zsigmondy primes}
\label{Section-A-Carlitz-module-analogue-of-Zsigmondy-primes}

In this section, we introduce a Carlitz module analogue of Zsigmondy primes for a pair of monic polynomials. We then prove a function field analogue of L\"uneburg's theorem which is crucial in the proofs of our main results. 

Using Goss \cite[Proposition 1.6.5 and Lemma 1.6.8]{Goss}, one obtains the following result.

\begin{proposition}
\label{Proposition-The-main-result-about-Carlitz-annihilator-of-u-and-wp}

Let $u$ be a nonzero polynomial in $A$, and let $\wp$ be a monic prime in $A$ such that $\wp$ does not divide $u$. Then there exists a unique monic polynomial $\fq$ of positive degree satisfying the following two conditions.
\begin{itemize}

\item [(CA1)] $C_{\fq}(u) \equiv 0 \pmod{\wp}$; and

\item [(CA2)] for any nonzero polynomial $m$ in $A$, $\fq$ divides $m$ if and only if $C_m(u) \equiv 0 \pmod{\wp}$.

\end{itemize}

\end{proposition}

In \cite{dongquan-JNT-145-2014}, I introduced the notion of the Carlitz annihilator of a monic prime to study congruences of primes dividing a Mersenne number in the function field setting. The following definition is a generalization of the notion of the Carlitz annihilator of a monic prime to a couple $(u, \wp)$, where $u$ is a nonzero polynomial and $\wp$ is a monic prime.

\begin{definition}
\label{Definition-The-Carlitz-annihilator-of-(u, wp)}

Let $u$ be a nonzero polynomial in $A$, and let $\wp$ be a monic prime in $A$. Let $P_{u, \wp}$ be the unique monic polynomial satisfying (CA1) and (CA2) in Proposition \ref{Proposition-The-main-result-about-Carlitz-annihilator-of-u-and-wp} if $\wp$ does not divide $u$, and let $P_{u, \wp} = 1$ if $\wp$ divides $u$. The monic polynomial $P_{u, \wp}$ is called \textit{the Carlitz annihilator of $(u, \wp)$}.

\end{definition}

For the rest of this paper, we always denote by $P_{u, \wp}$ the Carlitz annihilator of a couple $(u, \wp)$. The following result is immediate from Lemma \ref{Lemma-An-analogue-of-Fermat-little-theorem}, Proposition \ref{Proposition-The-main-result-about-Carlitz-annihilator-of-u-and-wp}, and Definition \ref{Definition-The-Carlitz-annihilator-of-(u, wp)}.

\begin{proposition}
\label{Proposition-The-1st-basic-results-about-the-Carlitz-annihilator-of-(u, wp)}

Let $u$ be a nonzero polynomial in $A$, and let $\wp$ be a monic prime in $A$. Let $P_{u, \wp}$ be the Carlitz annihilator of $(u, \wp)$. Then

\begin{itemize}

\item [(i)] $P_{u, \wp}$ divides $\wp - 1$;

\item [(ii)]  $P_{u, \wp} = 1$ if and only if $\wp$ divides $u$; 

\item [(iii)] for any nonzero polynomial $m$ in $A$, $P_{u, \wp}$ divides $m$  if and only if $C_m(u) \equiv 0 \pmod{\wp}$;

\item [(iv)] $C_n(u) \not\equiv 0 \pmod{\wp}$ for every nonzero polynomial $n \in A$ with $\deg(n) < \deg(P_{u, \wp})$.
\end{itemize}

\end{proposition}

\begin{proof}

Only part $(iv)$ needs a proof. If $P_{u, \wp} = 1$, then part $(iii)$ is trivial. Assume that $P_{u, \wp}$ is of positive degree, and take any nonzero polynomial $n \in A$ with $\deg(n) < \deg(P_{u, \wp})$. If $C_n(u) \equiv 0 \pmod{\wp}$, then we see from Proposition \ref{Proposition-The-main-result-about-Carlitz-annihilator-of-u-and-wp} that $P_{u, \wp}$ divides $n$. This implies that $\deg(P_{u, \wp}) \le \deg(n)$, which is a contradiction. Hence $C_n(u) \not\equiv 0 \pmod{\wp}$, and therefore our contention follows immediately.

\end{proof}

The analogy between the multiplicative group $\bG_m$ and the Carlitz module $C$ that is explained in Section \ref{Section-Introduction} (see Table \ref{Table-The-analogy-between-the-multiplicative-group-G_m-and-the-Carlitz-module-C}) motivates the following notion. 

\begin{definition}
\label{Definition-The-1st-definition-of-Zsigmondy-primes}

Let $m, u$ be nonzero polynomials in $A$. A \textit{Zsigmondy prime} for $(u, m)$ is a monic prime $\wp$ such that $C_m(u) \equiv 0 \pmod{\wp}$ and $C_n(u) \not\equiv 0 \pmod{\wp}$ for every nonzero polynomial $n \in A$ with $\deg(n) < \deg(m)$. 

If a monic prime $\wp$ is not a Zsigmondy prime for $(u, m)$, we say that $\wp$ is a \textit{non-Zsigmondy prime for $(u, m)$.} 

\end{definition}

\begin{remark}
\label{Remark-The-remark-after-the-1st-definition-of-Zsigmondy-primes}

\begin{itemize}

\item []

\item [(i)] Write $m = \epsilon m_0$, where $\epsilon \in \bF_q^{\times}$, and $m_0$ is a monic polynomial in $A$. Proposition \ref{Proposition-The-1st-basic-results-about-the-Carlitz-annihilator-of-(u, wp)}$(iv)$ implies that Definition \ref{Definition-The-1st-definition-of-Zsigmondy-primes} is equivalent to the following notion: \textit{A monic prime $\wp$ is a Zsigmondy prime for a pair $(u, m)$ if the Carlitz annihilator of $(u, \wp)$ is exactly $m_0$.} In particular, when $m$ is monic, this means that the Carlitz annihilator of $(u, \wp)$ is exactly $m$.

\item [(ii)] For any nonzero polynomials $m, u \in A$, write $m = \epsilon m_0$ and $u = \delta u_0$, where $\epsilon, \delta \in \bF_q^{\times}$, and $m_0, u_0$ are monic polynomials. We see that $C_m(u) = \epsilon\delta C_{m_0}(u_0)$, and $C_n(u) = \delta C_n(u_0)$ for any polynomial $n \in A$. Thus a prime $\wp$ is a Zsigmondy prime for $(u, m)$ if and only if it is a Zsigmondy prime for $(u_0, m_0)$.

\item [(iii)] Assume that $m, u$ are nonzero polynomials in $A$ such that $\deg(m) = \deg(u) = 0$. Then $m, u$ are units in $\bF_q^{\times}$. Hence $C_m(u) = mu \in \bF_q^{\times}$, which proves that there exist no Zsigmondy primes for $(u, m)$.

\end{itemize}

\end{remark}

By Remark \ref{Remark-The-remark-after-the-1st-definition-of-Zsigmondy-primes}, it suffices to study Zsigmondy primes for pairs $(u, m)$, where $m, u$ are monic polynomials such that at least one of them is of positive degree. For the rest of this paper, whenever we study Zsigmondy primes for a pair $(u, m)$, we will always assume that $m, u$ satisfy such conditions.

\begin{remark}
\label{Remark-A-Zsigmondy-prime-for-(u,m)-divides-Phi-m-(x)}

\begin{itemize}

\item []

\item [(i)] If $\wp$ is a Zsigmondy prime for $(u, m)$, then $\wp$ divides $C_m(u) = \prod_{\substack{n \mid m \\ \text{$n$ monic}} }\Psi_n(u)$. If $\wp$ divides $\Psi_n(u)$ for some monic polynomial $n$ dividing $m$ with $n \ne m$, then $C_n(u) = \prod_{\substack{b \mid n \\ \text{$b$ monic}} }\Psi_b(u) \equiv 0 \pmod{\wp}$. Note that $\deg(n) < \deg(m)$, and thus the last congruence implies that $\wp$ is not a Zsigmondy prime for $(u, m)$, which is a contradiction. Hence $\gcd(\wp, \Psi_n(u)) = 1$ for any monic polynomial $n$ dividing $m$ with $n \ne m$. Therefore $\wp$ divides $\Psi_m(u)$. 

\item [(ii)] The above remark also implies that if $\wp$ does not divide $\Psi_m(u)$, then $\wp$ is a non-Zsigmondy prime for $(u, m)$. The next result gives a simple criterion for testing whether a monic prime $\wp$ is a non-Zsigmondy prime for a given pair $(u, m)$. This is a function field analogue of L\"uneburg's theorem (see L\"uneburg \cite[Satz 1]{Luneburg} or Roitman \cite[Proposition 2]{Roitman}).

\end{itemize}

\end{remark}

\begin{theorem}
\label{Theorem-The-1st-version-of-Zsigmondy-Theorem}

Let $m, u$ be monic polynomials in $A$ such that $m$ is of positive degree. Let $\wp$ be a monic prime dividing $\Psi_m(u)$, and let $P_{u, \wp}$ be the Carlitz annihilator of $(u, \wp)$. Assume that condition $(D*)$ in Corollary \ref{Corollary-C_m(u)-is-nonzero} is satisfied. Then
\begin{itemize}

\item [(i)] $\wp$ is a non-Zsigmondy prime for $(u, m)$ if and only if $\wp$ divides $m$.

\item [(ii)] if $\wp$ is a non-Zsigmondy prime for $(u, m)$, then $m = P_{u, \wp} \wp^s$ for some positive integer $s$. Furthermore $\wp^2$ does not divide $\Psi_m(u)$ unless $q = 2$ and $\deg(\wp) = 1$.

\end{itemize}

\end{theorem}

\begin{proof}

Since $\wp$ divides $\Psi_m(u)$, it follows from Proposition \ref{Proposition-The-2nd-basic-result-of-cyclotomic-polynomials} that $\wp$ divides $C_m(u)$. From Proposition \ref{Proposition-The-1st-basic-results-about-the-Carlitz-annihilator-of-(u, wp)}$(iv)$, one sees that $P_{u, \wp}$ divides $m$. In particular this implies that $\deg(P_{u, \wp}) < \deg(m)$.

We now prove part $(i)$. We first show that the ``only if'' part of $(i)$ is true. Indeed assume that $\wp$ is a non-Zsigmondy prime for $(u, m)$. We know from Remark \ref{Remark-The-remark-after-the-1st-definition-of-Zsigmondy-primes} that $m \ne P_{u, \wp}$. Write
\begin{align}
\label{Equation-The-equation-of-m-in-terms-of-P-of-u-wp-and-n}
m = P_{u, \wp}^r n,
\end{align}
where $r \in \bZ_{>0}$, and $n$ is a monic element in $A$ with $\gcd(P_{u, \wp}, n) = 1$. Since $\deg(P_{u, \wp}) < \deg(m)$, it follows that $r \ge 2$ or $\deg(n) \ge 1$.

Let $l$ be a monic prime of positive degree such that $l$ divides $m$ and $P_{u, \wp}$ divides $m/l$. Then Proposition \ref{Proposition-The-1st-basic-results-about-the-Carlitz-annihilator-of-(u, wp)}$(iii)$ implies that
\begin{align}
\label{Equation-The-1st-equation-in-Zsigmondy-theorem}
C_{m/l}(u) \equiv 0 \pmod{\wp}.
\end{align}
Lemma \ref{Corollary-The-4th-basic-result-of-cyclotomic-polynomials} tells us that $\Psi_m(x)$ divides $C_m(x)/C_{m/l}(x)$ in $A[x]$. By Corollary \ref{Corollary-C_m(u)-is-nonzero}, $C_{m/l}(u) \ne 0$, and thus $\Psi_m(u)$ divides $C_m(u)/C_{m/l}(u)$. Since $\wp$ divides $\Psi_m(u)$, we deduce that
\begin{align}
\label{Equation-The-2nd-equation-in-Zsigmondy-theorem}
\dfrac{C_m(u)}{C_{m/l}(u)} \equiv 0 \pmod{\wp}.
\end{align}

By \cite[Proposition {\bf12.11}]{Rosen}, one can write $C_l(x) \in A[x]$ in the form
\begin{align*}
C_l(x) = lx + [l, 1]x^q + \ldots + [l, \deg(l) - 1]x^{q^{\deg(l) - 1}} + x^{q^{\deg(l)}},
\end{align*}
where $[l, i]$ is a polynomial of degree $q^i(\deg(l) - i)$ for each $1 \le i \le \deg(l) - 1$. By $(\ref{Equation-The-1st-equation-in-Zsigmondy-theorem})$, we deduce that
\begin{align*}
\dfrac{C_m(u)}{C_{m/l}(u)} &= \dfrac{C_l(C_{m/l}(u))}{C_{m/l}(u)} \\
&= l + [l, 1](C_{m/l}(u))^{q - 1} + \ldots  + (C_{m/l}(u))^{q^{\deg(l)} - 1} \equiv l \pmod{\wp},
\end{align*}
and therefore it follows from $(\ref{Equation-The-2nd-equation-in-Zsigmondy-theorem})$ that $l \equiv 0 \pmod{\wp}$. Since $l, \wp$ are monic primes, we deduce that $l = \wp$.

\textit{In summary, we have proved that if $l$ is a monic prime of positive degree such that $l$ divides $m$, and $P_{u, \wp}$ divides $m/l$, then $l = \wp$.}

We now prove that $r = 1$ and $n = \wp^s$ for some positive integer $s$ in the equation (\ref{Equation-The-equation-of-m-in-terms-of-P-of-u-wp-and-n}) of $m$. Indeed if $r \ge 2$, then letting $l = P_{u, \wp}$, we see that $\dfrac{m}{l} = \dfrac{m}{P_{u, \wp}} = P_{u, \wp}^{r - 1}n \equiv 0 \pmod{P_{u, \wp}}$. Thus repeating the same arguments as before, we deduce that $P_{u, \wp} = l \equiv 0 \pmod{\wp}$. On the other hand, Proposition \ref{Proposition-The-1st-basic-results-about-the-Carlitz-annihilator-of-(u, wp)}$(i)$ tells us that $P_{u, \wp}$ divides $\wp - 1$, and hence $P_{u, \wp} = \wp = \wp - 1$, which is a contradiction. Thus $r = 1$.

Now take any monic prime $l$ dividing $n$. From $(\ref{Equation-The-equation-of-m-in-terms-of-P-of-u-wp-and-n})$, one sees that $\dfrac{m}{l} = P_{u, \wp}(n/l) \equiv 0 \pmod{P_{u, \wp}}$. The same arguments as before yields that $l = \wp$, and thus $n = \wp^s$ for some $s \in \bZ_{>0}$. Therefore
\begin{align}
\label{Equation-The-2nd-representation-of-m-in-Zsigmondy-theorem}
m = P_{u, \wp}\wp^{s},
\end{align}
which implies that $\wp$ divides $m$. Hence the ``only if'' part of $(i)$ follows.

Suppose now that $\wp$ divides $m$. Assume the contrary, i.e., $\wp$ is a Zsigmondy prime for $(u, m)$. Then Remark \ref{Remark-The-remark-after-the-1st-definition-of-Zsigmondy-primes}$(i)$ tells us that $m = P_{u, \wp}$. Since $\wp$ divides $m$, we deduce that $\deg(\wp) \le \deg(m) = \deg(P_{u, \wp})$. On the other hand, we know from Proposition \ref{Proposition-The-1st-basic-results-about-the-Carlitz-annihilator-of-(u, wp)}$(i)$ that $P_{u, \wp}$ divides $\wp - 1$, and thus $\deg(P_{u, \wp}) \le \deg(\wp - 1) = \deg(\wp)$.
Therefore $\deg(P_{u, \wp}) = \deg(\wp) = \deg(\wp - 1)$, and hence $m = P_{u, \wp} = \wp - 1$, which is a contradiction since $\wp$ divides $m$. Thus $\wp$ is a non-Zsigmondy prime for $(u, m)$.

We now prove part $(ii)$. In the proof of the ``only if'' part of $(i)$ above, we have showed that if $\wp$ is a non-Zsigmondy prime for $(u, m)$, then $m$ is of the form (\ref{Equation-The-2nd-representation-of-m-in-Zsigmondy-theorem}), which proves the first part of $(ii)$. 

For the last part of $(ii)$, using \cite[Proposition {\bf12.11}]{Rosen}, one can write $C_{\wp}(x) \in A[x]$ in the form
\begin{align}
\label{Eqn-C-wp(x)-in-the-Luneburg-theorem-in-the-proof-of-part-(ii)}
C_{\wp}(x) = \wp x + [\wp, 1]x^q + \ldots + [\wp, d - 1]x^{q^{d - 1}} + x^{q^d},
\end{align}
where $d = \deg(\wp)$, and $[\wp, i]$ is a polynomial of degree $q^i(d - i)$ in $A$ for each $1 \le i \le d - 1$. It is well-known \cite[Proposition {\bf2.4}]{Hayes} that $C_{\wp}(x)$ is an Eisenstein polynomial, that is, $[\wp, i]$ is divisible by $\wp$ for each $1 \le i \le d - 1$.

By the first part of $(ii)$ and (\ref{Equation-The-2nd-representation-of-m-in-Zsigmondy-theorem}), we know that $P_{u, \wp}$ divides the polynomial $m/\wp$, and it thus follows from Proposition \ref{Proposition-The-1st-basic-results-about-the-Carlitz-annihilator-of-(u, wp)}$(iii)$ that $C_{m/\wp}(u) \equiv 0 \pmod{\wp}$.

If $d = \deg(\wp) > 1$, then we know that $q^{i} - 1 \ge 1$ for every $1 \le i \le d - 1$, and thus
\begin{align*}
[\wp, i]C_{m/\wp}(u)^{q^i - 1} \equiv 0 \pmod{\wp^2}
\end{align*}
for every $1 \le i \le  d - 1$. Furthermore since $q^d - 1 \ge q^2 - 1 \ge 3$, we deduce that $(C_{m/\wp}(u))^{q^d - 1} \equiv 0 \pmod{\wp^2}$. By Corollary \ref{Corollary-C_m(u)-is-nonzero}, one knows that $C_{m/\wp}(u) \ne 0$, and therefore we deduce from (\ref{Eqn-C-wp(x)-in-the-Luneburg-theorem-in-the-proof-of-part-(ii)}) that
\begin{align*}
\dfrac{C_m(u)}{C_{m/\wp}(u)} &= \dfrac{C_{\wp}(C_{m/\wp}(u))}{C_{m/\wp}(u)} \\
&= \wp  + [\wp, 1](C_{m/\wp}(u))^{q - 1} + \ldots  + (C_{m/\wp}(u))^{q^d - 1} \equiv \wp \pmod{\wp^2}.
\end{align*}

Since $C_{m/\wp}(u) \ne 0$, Lemma \ref{Corollary-The-4th-basic-result-of-cyclotomic-polynomials} tells us that $\Psi_m(u)$ divides $C_m(u)/C_{m/\wp}(u)$. If $\wp^2$ divides $\Psi_m(u)$, then $\wp \equiv \dfrac{C_m(u)}{C_{m/\wp}(u)} \equiv 0 \pmod{\wp^2}$, which is a contradiction. Therefore $\wp^2$ does not divide $\Psi_m(u)$.

If $q > 2$, then $q^i - 1 \ge 2$ for every $1 \le i \le  d$, and thus $(C_{m/\wp}(u))^{q^i - 1} \equiv 0 \pmod{\wp^2}$ for every $1 \le i \le d$. Repeating in the same arguments as above, we also obtain that $\wp^2$ does not divide $\Psi_m(u)$. 

\end{proof}

When $q > 2$, we see that condition $(D*)$ in Corollary \ref{Corollary-C_m(u)-is-nonzero} is trivially satisfied. For the rest of this paper, we will always assume that $q > 2$, and thus it is worth restating Theorem \ref{Theorem-The-1st-version-of-Zsigmondy-Theorem} with the assumption $q > 2$ in place of $(D*)$.

\begin{corollary}
\label{C-Luneburg-theorem-when-q->-2}

Assume that $q > 2$. Let $m, u$ be monic polynomials in $A$ such that $m$ is of positive degree. Let $\wp$ be a monic prime dividing $\Psi_m(u)$, and let $P_{u, \wp}$ be the Carlitz annihilator of $(u, \wp)$. Then
\begin{itemize}

\item [(i)] $\wp$ is a non-Zsigmondy prime for $(u, m)$ if and only if $\wp$ divides $m$.

\item [(ii)] if $\wp$ is a non-Zsigmondy prime for $(u, m)$, then $m = P_{u, \wp} \wp^s$ for some positive integer $s$. Furthermore $\wp^2$ does not divide $\Psi_m(u)$.

\end{itemize}

\end{corollary}

\section{Bang--Zsigmondy's theorem in characteristic $p \ne 2$}
\label{Section-Zsigmondy-theorem-in-characteristic-p-not-2}

In this section, we prove a function field analogue of Bang--Zsigmondy's theorem in the case when $p \ne 2$ (see Theorem \ref{Theorem-The-full-general-version-of-Zsigmondy-Theorem-with-p-not--equal-to-2} below). We begin by proving an elementary but very useful result that will play a key role in the proofs of Theorem \ref{Theorem-The-full-general-version-of-Zsigmondy-Theorem-with-p-not--equal-to-2} and Theorem \ref{Theorem-Zsigmondy-theorem-in-characteristic-two}.

\begin{lemma}
\label{Lemma-The-lower-bound-of-Phi(m)}

Assume that $q > 2$. Let $m$ be a polynomial in $A$ of positive degree. Then $\Phi(m) \ge (q - 1)\deg(m)$, where $\Phi(\cdot)$ denotes the function field analogue of the Euler $\phi$-function (see Subsection \ref{Notation} for its definition).

\end{lemma}

\begin{proof}

Since $\Phi(m) = \Phi(\epsilon m)$ for any $\epsilon \in \bF_q^{\times}$, one can, without loss of generality, assume that $m$ is a monic polynomial. We first prove Lemma \ref{Lemma-The-lower-bound-of-Phi(m)} when $m = P^s$, where $P$ is a monic prime and $s \in \bZ_{>0}$. Indeed, set $d = \deg(P) \ge 1$. If $s = 1$, using Bernoulli's inequality (see \cite[Theorem 42]{Hardy-Littlewood-Polya}), we deduce that
\begin{align*}
\Phi(m) = \Phi(P) = q^d - 1 = (1 + (q - 1))^d - 1 \ge 1 + (q - 1)d - 1 = (q - 1)d = (q - 1)\deg(m).
\end{align*}

If $s > 1$, then $q^{d(s - 1)} \ge 1 + (q - 1)d(s - 1)$, and thus $q^{d(s - 1)} > 1 + (s - 1) = s$. Therefore
\begin{align*}
\Phi(m) = \Phi(P^s) &= q^{\deg(P^s)} - q^{\deg(P^{s - 1})} = q^{d(s - 1)}(q^d - 1) \\
&> s((q - 1)d) = (q - 1)ds = (q - 1)\deg(P^s) = (q - 1)\deg(m),
\end{align*}
which proves Lemma \ref{Lemma-The-lower-bound-of-Phi(m)} for $m = P^s$.

Now let $m$ be any monic polynomial of positive degree. Take a monic prime $P$ dividing $m$, and write $m = nP^s$ for some $s \in \bZ_{>0}$, where $n$ is a monic polynomial with $\gcd(n, P) = 1$. If $n = 1$, then Lemma \ref{Lemma-The-lower-bound-of-Phi(m)} is true as shown above. If $\deg(n) \ge 1$, the induction implies that $\Phi(n) \ge (q - 1)\deg(n) \ge 2$. Therefore $\Phi(P^s)(\Phi(n) - 1) \ge (q-1)\deg(P^s)(\Phi(n) - 1) \ge 2(\Phi(n) - 1) \ge \Phi(n)$, and thus $\Phi(m) = \Phi(n)\Phi(P^s) \ge \Phi(n) + \Phi(P^s)$. Hence
\begin{align*}
\Phi(m) \ge \Phi(n) + \Phi(P^s) \ge (q - 1)\deg(n) + (q - 1)\deg(P^s) = (q - 1)\deg(nP^s) = (q- 1)\deg(m).
\end{align*}

\end{proof}

\begin{lemma}
\label{L-nonvanishing-of-Psi_m(u)}

Assume that $q > 2$. Let $m, u$ be monic polynomials such that $m$ is of positive degree. Then $\Psi_m(u)$ is of positive degree.

\end{lemma}

\begin{proof}

If $u$ is of positive degree, combining Lemmas \ref{Lemma-The-degrees-of-cyclotomic-polynomials} and \ref{Lemma-The-lower-bound-of-Phi(m)} yields 
\begin{align*}
\deg(\Psi_m(u)) = \deg(u)\Phi(m) \ge (q - 1)\deg(m) \ge q - 1 \ge 2,
\end{align*}
and thus $\Psi_m(u)$ is of positive degree. If $u$ is a unit in $\bF_q^{\times}$, using Lemma \ref{Lemma-The-degrees-of-cyclotomic-polynomials-when-u-is-a-unit} and repeating the same arguments as above, one also obtains that $\Psi_m(u)$ is of positive degree.

\end{proof}

The next result is crucial in the proof of Theorem \ref{Theorem-The-full-general-version-of-Zsigmondy-Theorem-with-p-not--equal-to-2}.

\begin{corollary}
\label{Corollary-The-description-of-non-Zsigmondy-primes}

Assume that $p \ne 2$ (recall that $p$ is the characteristic of $\bF_q$). Let $m, u$ be monic polynomials in $A$ such that $m$ is of positive degree. Assume that there are no Zsigmondy primes for $(u, m)$. Then $\Psi_m(u) = \epsilon\wp$ for some unit $\epsilon \in \bF_q^{\times}$ and some monic prime $\wp$.

\end{corollary}

\begin{proof}

Note that $q > 2$. By Lemma \ref{L-nonvanishing-of-Psi_m(u)}, $\Psi_m(u)$ is of positive degree, and hence there is a monic prime $\wp$ dividing $\Psi_m(u)$. Since $\wp$ is a non-Zsigmondy prime for $(u, m)$, applying Corollary \ref{C-Luneburg-theorem-when-q->-2}$(ii)$, one sees that $m = P_{u, \wp}\wp^s$ for some $s \in \bZ_{>0}$, where $P_{u, \wp}$ is the Carlitz annihilator of $(u, \wp)$.

We now prove that $\wp$ is the only monic prime factor of $\Psi_m(u)$. Indeed, assume the contrary, that is, there exists another monic prime $P$ dividing $\Psi_m(u)$ with $P \ne \wp$. By assumption, $P$ is a non-Zsigmondy prime for $(u, m)$. Hence it follows from Corollary \ref{C-Luneburg-theorem-when-q->-2}$(i)$ that $P$ divides $m$. Since $\gcd(P, \wp) = 1$ and $m = P_{u, \wp}\wp^s$, we see that $P$ divides $P_{u, \wp}$. Therefore by Proposition \ref{Proposition-The-1st-basic-results-about-the-Carlitz-annihilator-of-(u, wp)}$(i)$, we deduce that
\begin{align}
\label{Equation-The-1st-equation-in-the-corollary-about-non-Zsigmondy-primes}
\wp - 1 \equiv 0 \pmod{P}.
\end{align}
In particular, this implies that $\deg(P) \le \deg(\wp - 1) = \deg(\wp)$. Exchanging the roles of $\wp$ and $P$ and using the same arguments as above, we have
\begin{align}
\label{Equation-The-2nd-equation-in-the-corollary-about-non-Zsigmondy-primes}
P - 1 \equiv 0 \pmod{\wp}.
\end{align}
Furthermore this implies that $\deg(\wp) \le \deg(P - 1) = \deg(P)$, and thus $\deg(P) = \deg(\wp)$. Therefore we deduce from $(\ref{Equation-The-1st-equation-in-the-corollary-about-non-Zsigmondy-primes})$ and $(\ref{Equation-The-2nd-equation-in-the-corollary-about-non-Zsigmondy-primes})$ that $P = \wp - 1$ and $P - 1 = \wp$. Hence $P = \wp - 1 = P - 2$, and thus $2 = 0$, which is a contradiction since $p \ne 2$. Thus $\wp$ is the only monic prime factor of $\Psi_m(u)$, and hence $\Psi_m(u) = \epsilon \wp^e$ for some $\epsilon \in \bF_q^{\times}$ and some $e \in \bZ_{>0}$. By Corollary \ref{C-Luneburg-theorem-when-q->-2}$(ii)$, $\wp^2$ does not divide $\Psi_m(u)$, and hence $e = 1$, which proves Corollary \ref{Corollary-The-description-of-non-Zsigmondy-primes}.

\end{proof}

The following result is a Carlitz module analogue of Bang--Zsigmondy's theorem in characteristic $p \ne 2$. An analogue of Bang--Zsigmondy's theorem in characteristic two will be proved in Section \ref{Section-Zsigmondy-theorem-in-characteristic-two}.

\begin{theorem}
\label{Theorem-The-full-general-version-of-Zsigmondy-Theorem-with-p-not--equal-to-2}

Assume that $p \ne 2$. Let $m, u$ be monic polynomials in $A$ such that at least one of them is of positive degree. Then there exists a Zsigmondy prime for $(u, m)$ except exactly in the following case:
\begin{itemize}

\item [(EC1)] $q = 3$, $u = 1$, and $m = (\wp - 1)\wp$, where $\wp$ is an arbitrary monic prime of degree 1 in $\bF_3[T]$.

\end{itemize}

\end{theorem}

The proof of Theorem \ref{Theorem-The-full-general-version-of-Zsigmondy-Theorem-with-p-not--equal-to-2} will follow immediately from the next two lemmas.

\begin{lemma}
\label{Lemma-The-1st-lemma-for-the-full-Zsigmondy-Theorem-with-p-not-equal-to-2}

Assume that $p \ne 2$. Let $m, u$ be monic polynomials in $A$ such that $u$ is of positive degree. Then there exists a Zsigmondy prime for $(u, m)$.

\end{lemma}

\begin{proof}

Note that $q > 2$. If $m = 1$, then $\wp$ is a Zsigmondy prime for $(u, m)$ for each monic prime $\wp$ dividing $u$. Suppose now that $\deg(m) > 0$. Assume the contrary, i.e, there exist no Zsigmondy primes for $(u, m)$. By Corollary \ref{Corollary-The-description-of-non-Zsigmondy-primes}, $\Psi_m(u) = \epsilon \wp$ for some unit $\epsilon \in \bF_q^{\times}$ and some monic prime $\wp$.

Since $u$ is of positive degree, one sees from Lemma \ref{Lemma-The-degrees-of-cyclotomic-polynomials} that $\deg(\wp) = \deg(\Psi_m(u)) = \deg(u)\Phi(m)$. Since $\wp$ is a non-Zsigmondy prime for $(u, m)$, we deduce from Corollary \ref{C-Luneburg-theorem-when-q->-2}$(i)$ that $\wp$ divides $m$, and thus $\deg(\wp) \le \deg(m)$. Therefore $\deg(u)\Phi(m) \le \deg(m)$.

On the other hand, Lemma \ref{Lemma-The-lower-bound-of-Phi(m)} tells us that $\Phi(m) \ge (q - 1)\deg(m)$, and therefore
 \begin{align*}
 \deg(m) \ge \deg(u)\Phi(m) \ge \Phi(m) \ge (q - 1)\deg(m) \ge 2\deg(m).
 \end{align*}
 Hence $\deg(m) = 0$, which is a contradiction. Thus there exists a Zsigmondy prime for $(u, m)$.

 \end{proof}

 \begin{lemma}
\label{Lemma-The-2nd-lemma-for-the-full-Zsigmondy-Theorem-with-p-not-equal-to-2}

Assume that $p \ne 2$. Let $m$ be a monic polynomial in $A$ such that $m$ is of positive degree. Then there exists a Zsigmondy prime for $(1, m)$ except exactly in the following case:
\begin{itemize}

\item [(EC1)] $q = 3$, and $m = (\wp - 1)\wp$, where $\wp$ is an arbitrary monic prime of degree 1 in $\bF_3[T]$.

\end{itemize}

\end{lemma}

\begin{proof}

Assume first that $q = 3$, and $m = (\wp - 1)\wp$, where $\wp$ is an arbitrary monic prime of degree 1 in $\bF_3[T]$. Table \ref{Table-Nonexistence-of-Zsigmondy-primes-for-(1,m)-when-q-is-3} tells us that there exist no Zsigmondy primes for $(1, m)$.

\begin{table}
\caption{Nonexistence of Zsigmondy primes for $(1,m)$, where $m = (\wp - 1)\wp$ for some monic prime $\wp$ of degree 1 in $\bF_3[T]$.}
\label{Table-Nonexistence-of-Zsigmondy-primes-for-(1,m)-when-q-is-3}

\begin{center}
\begin{tabular}{|l|c|p{5cm}|p{5cm}|}

\hline

$\wp$ & $m = (\wp - 1)\wp$ & The prime factorization of $C_m(1)$ & Nonexistence of Zsigmondy primes for $(1, m)$ \\

\hline

\hline

$T$ & $(T - 1)T$ & $C_m(1) = \wp_1^2\wp_2$, where $\wp_1 = T$ and $\wp_2 = T + 1$.  & There are no Zsigmondy primes for $(1, m)$ since $\wp_1 = C_{T - 1}(1)$ and  $C_T(1) = \wp_2$. \\

\hline

$T + 1$ & $T(T + 1)$ &  $C_m(1) = \wp_1^2\wp_2$, where $\wp_1 = T + 1$ and $\wp_2 = T + 2$. & There are no Zsigmondy primes for $(1, m)$ since $\wp_1 = C_{T}(1)$ and  $C_{T + 1}(1) = \wp_2$.\\

\hline

$T + 2$ & $(T + 1)(T + 2)$ &  $C_m(1) = \wp_1\wp_2^2$, where $\wp_1 = T$ and $\wp_2 = T + 2$. & There are no Zsigmondy primes for $(1, m)$ since $\wp_1 = C_{T - 1}(1)$ and  $C_{T + 1}(1) = \wp_2$\\

\hline

\end{tabular}
\end{center}
\end{table}

Suppose now that we are not in the exceptional case (EC1) in Lemma \ref{Lemma-The-2nd-lemma-for-the-full-Zsigmondy-Theorem-with-p-not-equal-to-2}, that is, either $q = 3$ and $m \ne (\wp - 1)\wp$ for any monic prime $\wp$ of degree 1 in $\bF_3[T]$ or $q \ne 3$. Assume the contrary, i.e., there exist no Zsigmondy primes for $(1, m)$. It then follows from Corollary \ref{Corollary-The-description-of-non-Zsigmondy-primes} that
\begin{align}
\label{Equation-The-1st-equation-in-the-full-Zsigmondy-theorem-with-p-not-equal-to-2}
 \Psi_m(1) = \epsilon \wp
 \end{align}
 for some $\epsilon \in \bF_q^{\times}$ and some monic prime $\wp$. Since $\wp$ is a non-Zsigmondy prime for $(1, m)$, we know from Corollary \ref{C-Luneburg-theorem-when-q->-2}$(ii)$ that
 \begin{align}
 \label{Equation-The-2nd-equation-in-the-full-Zsigmondy-theorem-with-p-not-equal-to-2}
 m = P_{1, \wp}\wp^s,
 \end{align}
 where $P_{1, \wp}$ is the Carlitz annihilator of $(1, \wp)$, and $s \in \bZ_{>0}$. By Proposition \ref{Proposition-The-1st-basic-results-about-the-Carlitz-annihilator-of-(u, wp)}$(ii)$, $\deg(P_{1, \wp}) > 0$.

By Lemma \ref{Lemma-The-degrees-of-cyclotomic-polynomials-when-u-is-a-unit}, we see that
\begin{align}
\label{Equation-The-3rd-equation-in-the-full-Zsigmondy-theorem-with-p-not-equal-to-2}
\deg(\Psi_m(1)) = \dfrac{\Phi(m) + \delta}{q},
\end{align}
where $\delta$ is an integer in $\{0, \pm 1\}$. Recall from Proposition \ref{Proposition-The-1st-basic-results-about-the-Carlitz-annihilator-of-(u, wp)}$(i)$ that $P_{1, \wp}$ divides $\wp - 1$, and thus $\gcd(P_{1, \wp}, \wp) = 1$. Note that $q > 2$, and thus Lemma \ref{Lemma-The-lower-bound-of-Phi(m)} yields
\begin{align*}
\Phi(m) = \Phi(P_{1, \wp})\Phi(\wp^s) \ge (q - 1)\deg(P_{1, \wp})(q - 1)\deg(\wp^s) \ge s(q - 1)^2\deg(P_{1, \wp})\deg(\wp).
\end{align*}
Since $\deg(P_{1, \wp}) \ge 1$ and $s \ge 1$, it follows from the above inequalities that
\begin{align}
\label{Inequality-The-1st-inequality-in-the-full-Zsigmondy-theorem-with-p-not-equal-to-2}
\Phi(m) \ge (q - 1)^2\deg(\wp).
\end{align}
Note that if equality in $(\ref{Inequality-The-1st-inequality-in-the-full-Zsigmondy-theorem-with-p-not-equal-to-2})$ occurs, then $\deg(P_{1, \wp}) = s = 1$.

Since $\delta \ge -1$, we deduce from $(\ref{Equation-The-3rd-equation-in-the-full-Zsigmondy-theorem-with-p-not-equal-to-2})$ and (\ref{Inequality-The-1st-inequality-in-the-full-Zsigmondy-theorem-with-p-not-equal-to-2}) that
\begin{align}
\label{Inequality-The-2nd-inequality-in-the-full-Zsigmondy-theorem-with-p-not-equal-to-2}
\deg(\Psi_m(1)) \ge \dfrac{(q - 1)^2\deg(\wp) - 1}{q}.
\end{align}
From $(\ref{Inequality-The-1st-inequality-in-the-full-Zsigmondy-theorem-with-p-not-equal-to-2})$ and the remark following $(\ref{Inequality-The-1st-inequality-in-the-full-Zsigmondy-theorem-with-p-not-equal-to-2})$, we see that if equality in $(\ref{Inequality-The-2nd-inequality-in-the-full-Zsigmondy-theorem-with-p-not-equal-to-2})$ occurs, then $\deg(P_{1, \wp}) = s = 1$.

By $(\ref{Equation-The-1st-equation-in-the-full-Zsigmondy-theorem-with-p-not-equal-to-2}), (\ref{Inequality-The-2nd-inequality-in-the-full-Zsigmondy-theorem-with-p-not-equal-to-2})$, we deduce that
\begin{align}
\label{Inequality-The-3rd-inequality-in-the-full-Zsigmondy-theorem-with-p-not-equal-to-2}
\dfrac{(q - 1)^2\deg(\wp) - 1}{q} \le \deg(\Psi_m(1)) = \deg(\wp),
\end{align}
and thus
\begin{align}
\label{Inequality-The-4th-inequality-in-the-full-Zsigmondy-theorem-with-p-not-equal-to-2}
((q-1)^2 - q)\deg(\wp) \le 1.
\end{align}
Note that equality in $(\ref{Inequality-The-4th-inequality-in-the-full-Zsigmondy-theorem-with-p-not-equal-to-2})$ occurs if and only if equality in $(\ref{Inequality-The-3rd-inequality-in-the-full-Zsigmondy-theorem-with-p-not-equal-to-2})$ occurs.

On the other hand, we see that $(q - 1)^2 - q \ge 1$ with equality if and only if $q = 3$. Since $\deg(\wp) \ge 1$, we deduce from $(\ref{Inequality-The-4th-inequality-in-the-full-Zsigmondy-theorem-with-p-not-equal-to-2})$ that $1 \le ((q-1)^2 - q)\deg(\wp) \le 1$, and therefore
\begin{align}
\label{Equation-The-4th-equation-in-the-full-Zsigmondy-theorem-with-p-not-equal-to-2}
((q-1)^2 - q)\deg(\wp) = 1.
\end{align}

Equation $(\ref{Equation-The-4th-equation-in-the-full-Zsigmondy-theorem-with-p-not-equal-to-2})$ implies that equality in (\ref{Inequality-The-4th-inequality-in-the-full-Zsigmondy-theorem-with-p-not-equal-to-2}) occurs, and thus $q = 3$ and $\deg(\wp) = 1$. Furthermore since equality in (\ref{Inequality-The-4th-inequality-in-the-full-Zsigmondy-theorem-with-p-not-equal-to-2}) occurs, equality in (\ref{Inequality-The-3rd-inequality-in-the-full-Zsigmondy-theorem-with-p-not-equal-to-2}) also occurs, which in turns implies that equality in (\ref{Inequality-The-2nd-inequality-in-the-full-Zsigmondy-theorem-with-p-not-equal-to-2}) occurs. Therefore by the remark following (\ref{Inequality-The-2nd-inequality-in-the-full-Zsigmondy-theorem-with-p-not-equal-to-2}), we see that $\deg(P_{1, \wp}) = s = 1$. Thus it follows from (\ref{Equation-The-2nd-equation-in-the-full-Zsigmondy-theorem-with-p-not-equal-to-2}) that $m = P_{1, \wp}\wp$, where $\deg(P_{1, \wp}) = \deg(\wp) = 1$. Since  $P_{1, \wp}, \wp - 1$ are monic primes and $P_{1, \wp}$ divides $\wp - 1$, we deduce that $P_{1, \wp} = \wp - 1$. Therefore $q = 3$ and $m = (\wp - 1)\wp$, where $\wp$ is a monic prime of degree $1$ in $\bF_3[T]$. This implies that we are in the exceptional case (EC1), which is a contradiction. Thus there exists a Zsigmondy prime for $(1, m)$.

\end{proof}

\section{Bang--Zsigmondy's theorem in characteristic two}
\label{Section-Zsigmondy-theorem-in-characteristic-two}

In this section, we will prove an analogue of Bang--Zsigmondy's theorem in characteristic two. Throughout this section, we assume that $p = 2$ and $q > 2$. 

\begin{lemma}
\label{Lemma-The-prime-factorization-of-Psi-m-u-when-there-are-no-Zsigmondy-primes-for-u-m-and-p=2}

Let $m, u$ be monic polynomials in $A$ such that $m$ is of positive degree. Assume that there are no Zsigmondy primes for $(u, m)$. Then either of the following is true:
\begin{itemize}

\item [(i)] $\Psi_m(u) = \epsilon\wp$ for some unit $\epsilon \in \bF_q^{\times}$ and some monic prime $\wp$;

\item [(ii)] $\Psi_m(u) = \epsilon\wp(\wp - 1)$, where $\epsilon \in \bF_q^{\times}$ and $\wp$ is a monic prime such that $\wp - 1$ is also a prime.

\end{itemize}
\end{lemma}

\begin{proof}

By Lemma \ref{L-nonvanishing-of-Psi_m(u)}, we know that $\Psi_m(u)$ is of positive degree. Hence there exists a monic prime $\wp$ dividing $\Psi_m(u)$. Since $\wp$ is a non-Zsigmondy prime for $(u, m)$, applying Corollary \ref{C-Luneburg-theorem-when-q->-2}$(ii)$, we deduce that
\begin{align}
\label{Equation-The-1st-equation-in-lemma-about-the-prime-factorization-of-Psi-m-u-with-p=2}
m = P_{u, \wp}\wp^s
\end{align}
for some positive integer $s$, where $P_{u, \wp}$ is the Carlitz annihilator of $(u, \wp)$.

If $\wp$ is the only prime factor of $\Psi_m(u)$, then $\Psi_m(u) = \epsilon \wp^e$ for some $\epsilon \in \bF_q^{\times}$ and $e \in \bZ_{>0}$. By Corollary \ref{C-Luneburg-theorem-when-q->-2}$(ii)$, $\wp^2$ does not divide $\Psi_m(u)$, and therefore $e = 1$. Thus $\Psi_m(u) = \epsilon \wp$.

If $\Psi_m(u)$ has at least two distinct prime factors, fix a monic prime factor $\wp$ of $\Psi_m(u)$, and take any monic prime $P$ dividing $\Psi_m(u)$ with $P \ne \wp$. We contend that $P = \wp - 1$. Indeed, by assumption, $P$ is a non-Zsigmondy prime for $(u, m)$. Hence Corollary \ref{C-Luneburg-theorem-when-q->-2}$(i)$ tells us that $P$ divides $m$. Since $\gcd(P, \wp) = 1$, we deduce from (\ref{Equation-The-1st-equation-in-lemma-about-the-prime-factorization-of-Psi-m-u-with-p=2}) that $P$ divides $P_{u, \wp}$. Therefore by Proposition \ref{Proposition-The-1st-basic-results-about-the-Carlitz-annihilator-of-(u, wp)}$(i)$, we have
\begin{align}
\label{Equation-The-1st-equation-in-the-lemma-about-the-prime-factorization-of-Psi-m-u-with-p=2}
\wp - 1 \equiv 0 \pmod{P}.
\end{align}
In particular, this implies that $\deg(P) \le \deg(\wp - 1) = \deg(\wp)$.

Exchanging the roles of $\wp$ and $P$ and using the same arguments as above, we deduce that
\begin{align}
\label{Equation-The-2nd-equation-in-the-lemma-about-the-prime-factorization-of-Psi-m-u-with-p=2}
P - 1 \equiv 0 \pmod{\wp}.
\end{align}
Furthermore this implies that $\deg(\wp) \le \deg(P - 1) = \deg(P)$, and it thus follows that $\deg(P) = \deg(\wp)$. Since $P, \wp$ are monic, we deduce from $(\ref{Equation-The-1st-equation-in-the-lemma-about-the-prime-factorization-of-Psi-m-u-with-p=2})$ that $P = \wp - 1$.

Thus we have proved that if $P$ is an arbitrary monic prime dividing $\Psi_m(u)$ such that $P \ne \wp$, then $P = \wp - 1$. In particular, this implies that $\wp - 1$ is prime, and $\wp, \wp - 1$ are the only monic prime factors of $\Psi_m(u)$. Therefore $\Psi_m(u)$ is of the form $\Psi_m(u) = \epsilon (\wp - 1)^r\wp^e$ for some $r, e \in \bZ_{>0}$ and $\epsilon \in \bF_q^{\times}$. Corollary \ref{C-Luneburg-theorem-when-q->-2}$(ii)$ now implies that $\wp^2$ does not divide $\Psi_m(u)$, and therefore $e = 1$. Similarly $(\wp - 1)^2$ does not divide $\Psi_m(u)$, and thus $r = 1$. Hence $\Psi_m(u) = \epsilon (\wp - 1)\wp$, which proves our contention.

\end{proof}

We now state an analogue of Bang--Zsigmondy's theorem in characteristic two.
\begin{theorem}
\label{Theorem-Zsigmondy-theorem-in-characteristic-two}

Let $m, u$ be monic polynomials in $A$ such that at least one of them is of positive degree. Then there exists a Zsigmondy prime for $(u, m)$ except exactly in the following case:
\begin{itemize}

\item [(EC2)] $q = 2^2$, $u = 1$, and $m = (\wp - 1)\wp$, where $\wp$ is an arbitrary monic prime of degree 1 in $\bF_{2^2}[T]$.

\end{itemize}
\end{theorem}

The proof of Theorem \ref{Theorem-Zsigmondy-theorem-in-characteristic-two} will follow immediately from the next two lemmas.

\begin{lemma}
\label{Lemma-The-1st-lemma-for-Zsigmondy-theorem-in-char-2}

Let $m, u$ be monic polynomials in $A$ such that $u$ is of positive degree. Then there exists a Zsigmondy prime for $(u, m)$.

\end{lemma}

\begin{proof}

If $m = 1$, then $\wp$ is a Zsigmondy prime for $(u, m)$ for each monic prime $\wp$ dividing $u$. Suppose now that $\deg(m) > 0$. Assume the contrary, i.e., there exist no Zsigmondy primes for $(u, m)$. By Lemma \ref{Lemma-The-prime-factorization-of-Psi-m-u-when-there-are-no-Zsigmondy-primes-for-u-m-and-p=2}, we know that either of the following is true:
\begin{itemize}

\item [(i)] $\Psi_m(u) = \epsilon \wp$ for some unit $\epsilon \in \bF_q^{\times}$ and some monic prime $\wp$;

\item [(ii)] $\Psi_m(u) = \epsilon (\wp - 1)\wp$ for some $\epsilon \in \bF_q^{\times}$, where both $\wp, \wp - 1$ are monic primes.
\end{itemize}

If $\Psi_m(u) = \epsilon \wp$ for some unit $\epsilon \in \bF_q^{\times}$ and some monic prime $\wp$, repeating the same arguments as in the proof of Lemma \ref{Lemma-The-1st-lemma-for-the-full-Zsigmondy-Theorem-with-p-not-equal-to-2}, we deduce that there exists a Zsigmondy prime for $(u, m)$.

Suppose now that $\Psi_m(u) = \epsilon (\wp - 1)\wp$ for some $\epsilon \in \bF_q^{\times}$, where both $\wp, \wp - 1$ are monic primes. Since $\wp$ divides $\Psi_m(u)$ and $\wp$ is a non-Zsigmondy prime for $(u, m)$, we deduce from Corollary \ref{C-Luneburg-theorem-when-q->-2}$(ii)$ that
\begin{align}
\label{Equation-The-1st-equation-in-the-1st-lemma-about-Zsigmondy-theorem-with-char-p=2}
m = P_{u, \wp}\wp^s,
\end{align}
where $P_{u, \wp}$ is the Carlitz annihilator of $(u, \wp)$ and $s \in \bZ_{>0}$. Similarly, Corollary \ref{C-Luneburg-theorem-when-q->-2}$(ii)$ yields
\begin{align}
\label{Equation-The-2nd-equation-in-the-1st-lemma-about-Zsigmondy-theorem-with-char-p=2}
m = P_{u, \wp - 1}(\wp - 1)^r,
\end{align}
where $P_{u, \wp - 1}$ is the Carlitz annihilator of $(u, \wp - 1)$ and $r \in \bZ_{>0}$. Since $(\wp - 1) - 1 = \wp$, we deduce from Proposition \ref{Proposition-The-1st-basic-results-about-the-Carlitz-annihilator-of-(u, wp)}$(i)$ that $P_{u, \wp - 1}$ divides $\wp$.

From (\ref{Equation-The-1st-equation-in-the-1st-lemma-about-Zsigmondy-theorem-with-char-p=2}), (\ref{Equation-The-2nd-equation-in-the-1st-lemma-about-Zsigmondy-theorem-with-char-p=2}), we get
\begin{align}
\label{Equation-The-3rd-equation-in-the-1st-lemma-about-Zsigmondy-theorem-with-char-p=2}
m = P_{u, \wp}\wp^s = P_{u, \wp - 1}(\wp - 1)^r.
\end{align}
By Proposition \ref{Proposition-The-1st-basic-results-about-the-Carlitz-annihilator-of-(u, wp)}$(i)$, $P_{u, \wp}$ divides $\wp - 1$, and hence $\deg(P_{u, \wp}) \le \deg(\wp - 1) = \deg(\wp)$. From (\ref{Equation-The-3rd-equation-in-the-1st-lemma-about-Zsigmondy-theorem-with-char-p=2}) and since $\gcd((\wp - 1)^r, \wp^s) = 1$, we deduce that $(\wp - 1)^r$ divides $P_{u, \wp}$. Hence
\begin{align*}
\deg(P_{u, \wp}) \le \deg(\wp) \le r\deg(\wp) = r\deg(\wp - 1) = \deg((\wp - 1)^r) \le \deg(P_{u, \wp}),
\end{align*}
which implies that $\deg(P_{u, \wp}) = \deg(\wp) = r\deg(\wp) = \deg((\wp - 1)^r)$. Therefore $r = 1$, and $(\wp - 1)^r = P_{u, \wp}$. Thus $P_{u, \wp} = \wp - 1$. Exchanging the roles of $\wp$ and $\wp - 1$, one can show that $s = 1$ and $P_{u, \wp - 1} = (\wp - 1) - 1 = \wp$, and thus $m = (\wp - 1)\wp$. Therefore $\Psi_m(u) = \epsilon (\wp - 1)\wp = \epsilon m$.

By Lemma \ref{Lemma-The-degrees-of-cyclotomic-polynomials}, $\deg(m) = \deg(\Psi_m(u)) = \deg(u)\Phi(m)$. Hence it follows from Lemma \ref{Lemma-The-lower-bound-of-Phi(m)} that
\begin{align*}
\deg(m) \ge \deg(u)(q - 1)\deg(m) \ge (q - 1)\deg(m) > \deg(m),
\end{align*}
which is a contradiction. Thus there exists a Zsigmondy prime for $(u, m)$.

\end{proof}

We now consider Theorem \ref{Theorem-Zsigmondy-theorem-in-characteristic-two} in the case when $m$ is of positive degree and $u = 1$.

\begin{lemma}
\label{Lemma-The-2nd-lemma-for-Zsigmondy-theorem-in-char-2}

Let $m$ be a monic polynomial in $A$ such that $m$ is of positive degree. Then there exists a Zsigmondy prime for $(1, m)$ except exactly in the following case:
\begin{itemize}

\item [(EC2)] $q = 2^2$, and $m = (\wp - 1)\wp$, where $\wp$ is an arbitrary monic prime of degree 1 in $\bF_{2^2}[T]$.

\end{itemize}

\end{lemma}

\begin{proof}

Assume that we are in the exceptional case (EC2), that is, $q = 2^2$ and $m = (\wp - 1)\wp$ for a monic prime $\wp$ of degree $1$ in $\bF_{2^2}[T]$. We see that $\bF_{2^2} = \bF_2(\omega)$, where $\omega$ is an element in the algebraic closure of $\bF_2$ such that $\omega^2 + \omega + 1 = 0$. We know that $\{T, T + 1, T + \omega, T + \omega^2\}$ consists of all monic primes of degree 1 in $\bF_{2^2}[T]$. Table \ref{Table-There-are-no-Zsigmondy-primes-for-(1,m)-when-q-equals-4} tells us that there are no Zsigmondy primes for $(1, m)$ in the exceptional case (EC2).

\begin{table}
\caption{Nonexistence of Zsigmondy primes for $(1,m)$, where $m = (\wp - 1)\wp$ for a monic prime $\wp$ of degree 1 in $\bF_{2^2}[T]$.}
\label{Table-There-are-no-Zsigmondy-primes-for-(1,m)-when-q-equals-4}
\begin{center}
\begin{tabular}{|l|c|p{5cm}|p{5cm}|}

\hline

$\wp$ & $m = (\wp - 1)\wp$ & \textbf{The prime factorization of $C_m(1)$} & \textbf{Nonexistence of Zsigmondy primes for $(1, m)$} \\

\hline

\hline

$T$ & $(T - 1)T$ & $C_m(1) = \wp_1^2\wp_2^2$, where $\wp_1 = T$ and $\wp_2 = T + 1$.  & There are no Zsigmondy primes for $(1, m)$ since $\wp_1 = C_{T - 1}(1)$ and  $\wp_2 = C_T(1)$. \\

\hline

$T + 1$ & $T(T + 1)$ &  $C_m(1) = \wp_1^2\wp_2^2$, where $\wp_1 = T$ and $\wp_2 = T + 1$.  & There are no Zsigmondy primes for $(1, m)$ since $\wp_1 = C_{T - 1}(1)$ and  $\wp_2 = C_T(1)$. \\

\hline

$T + \omega$ & $(T + \omega - 1)(T + \omega)$ &  $C_m(1) = \wp_1^2\wp_2^2$, where $\wp_1 = T + \omega$ and $\wp_2 = T + \omega^2$. & There are no Zsigmondy primes for $(1, m)$ since $\wp_1 = C_{T + \omega - 1}(1)$ and  $\wp_2 = C_{T + \omega^2 - 1}(1) $\\

\hline

$T + \omega^2$ & $(T + \omega^2 - 1)(T + \omega^2)$  &  $C_m(1) = \wp_1^2\wp_2^2$, where $\wp_1 = T + \omega$ and $\wp_2 = T + \omega^2$. & There are no Zsigmondy primes for $(1, m)$ since $\wp_1 = C_{T + \omega - 1}(1)$ and  $\wp_2 = C_{T + \omega^2 - 1}(1) $\\

\hline

\end{tabular}
\end{center}
\end{table}

Suppose now that we are not in the exceptional case (EC2) in Lemma \ref{Lemma-The-2nd-lemma-for-Zsigmondy-theorem-in-char-2}, that is, either $q = 2^2$ and $m \ne (\wp - 1)\wp$ for any monic prime $\wp$ of degree 1 in $\bF_{2^2}[T]$ or $q \ne 2^2$. Assume the contrary, i.e., there exist no Zsigmondy primes for $(1, m)$. By Lemma \ref{Lemma-The-prime-factorization-of-Psi-m-u-when-there-are-no-Zsigmondy-primes-for-u-m-and-p=2}, we know that either of the following is true:
\begin{itemize}

\item [(i)] $\Psi_m(1) = \epsilon \wp$ for some unit $\epsilon \in \bF_q^{\times}$ and some monic prime $\wp$;

\item [(ii)] $\Psi_m(1) = \epsilon (\wp - 1)\wp$ for some unit $\epsilon \in \bF_q^{\times}$, where both $\wp, \wp - 1$ are monic primes.

\end{itemize}
Repeating the same arguments as in the proof of Lemma \ref{Lemma-The-2nd-lemma-for-the-full-Zsigmondy-Theorem-with-p-not-equal-to-2}, we deduce that there exists a Zsigmondy prime for $(1, m)$ if $\Psi_m(1) = \epsilon \wp$ for some unit $\epsilon \in \bF_q^{\times}$ and some monic prime $\wp$.

Suppose now that $\Psi_m(1) = \epsilon (\wp - 1)\wp$ for some unit $\epsilon \in \bF_q^{\times}$, where both $\wp, \wp - 1$ are monic primes. Using the same arguments as in Lemma \ref{Lemma-The-1st-lemma-for-Zsigmondy-theorem-in-char-2}, we deduce that $m = (\wp - 1)\wp$ and hence $\Psi_m(1) = \epsilon m$. By Lemma \ref{Lemma-The-degrees-of-cyclotomic-polynomials-when-u-is-a-unit}, we see that $\deg(\Psi_m(1)) = \dfrac{\Phi(m) + \delta}{q}$ for some integer $\delta \in \{0, \pm 1\}$. Since $\wp -1, \wp$ are primes with $\gcd(\wp - 1, \wp) = 1$, Lemma \ref{Lemma-The-lower-bound-of-Phi(m)} yields
\begin{align*}
\Phi(m) = \Phi((\wp - 1)\wp) = \Phi(\wp - 1)\Phi(\wp) \ge (q - 1)\deg(\wp - 1)(q- 1)\deg(\wp) = (q-1)^2\deg(\wp)^2.
\end{align*}

Since $\delta \ge -1$ and $\deg(\wp) \ge 1$, we get
\begin{align*}
2\deg(\wp) = \deg(m) = \deg(\Psi_m(1)) \ge \dfrac{(q-1)^2\deg(\wp)^2 - 1}{q} \ge \dfrac{(q-1)^2\deg(\wp) - 1}{q},
\end{align*}
and therefore $\deg(\wp)((q-1)^2 - 2q) \le 1$. By assumption, $q \ge 4$, and thus $(q - 1)^2 - 2q = q^2 - 4q + 1 \ge 1$ with equality if and only if $q = 2^2$. Since $\deg(\wp) \ge 1$, we deduce that $1 \le \deg(\wp)((q-1)^2 - 2q) \le 1$, and thus $\deg(\wp)((q-1)^2 - 2q) = 1$. The last equation yields $\deg(\wp) = 1$ and $q = 2^2$. Hence $q = 2^2$, and $m = (\wp - 1)\wp$ for some monic prime $\wp$ of degree 1 in $\bF_{2^2}[T]$. Thus we are in the exceptional case (EC2), which is a contradiction. Therefore there exists a Zsigmondy prime for $(1, m)$.

\end{proof}

\section{An analogue of large Zsigmondy primes, and Feit's theorem in positive characteristic}
\label{Section-Large-Zsigmondy-primes}

In this section, we introduce a notion of large Zsigmondy primes in the function field context, and prove a function field analogue of Feit's theorem (see Theorem \ref{Theorem-The-main-theorem-about-large-Zsigmondy-primes}) which assures the existence of a large Zsigmondy prime for a pair $(u, m)$ of monic polynomials except some exceptional cases. We begin by recalling the notion of large Zsigmondy primes that was already mentioned in the introduction. 

\begin{definition}
\label{D-lagre-Zsigmondy-primes}

Let $m, u$ be monic polynomials in $A = \bF_q[T]$. A monic prime  $\wp$ is called a \textit{large Zsigmondy prime for $(u, m)$} if $\wp$ is a Zsigmondy prime for $(u, m)$, and either $\deg(\wp) > \deg(m)$ or $\wp^2$ divides $C_m(u)$.

\end{definition}

\begin{remark}
\label{R-large-Zsigmondy-primes}

Let $\wp$ be a Zsigmondy prime for $(u, m)$. We know from Remark \ref{Remark-The-remark-after-the-1st-definition-of-Zsigmondy-primes}$(i)$ that $m$ is the Carlitz annihilator of $(u, \wp)$. By Proposition \ref{Proposition-The-1st-basic-results-about-the-Carlitz-annihilator-of-(u, wp)}$(i)$, we deduce that $\deg(m) \le \deg(\wp - 1) = \deg(\wp)$. Thus $\wp$ is not a large Zsigmondy prime for $(u, m)$ if and only if the following are true:
\begin{itemize}

\item [(i)] $\deg(\wp) = \deg(m)$; and

\item [(ii)] $\wp^2$ does not divide $C_m(u)$.

\end{itemize}

\end{remark}

The next result plays a central role in the proof of Theorem \ref{Theorem-The-main-theorem-about-large-Zsigmondy-primes}.

\begin{corollary}
\label{Corollary-Description-of-large-Zsigmondy-primes}

Assume that $q > 2$. Let $m, u$ be monic polynomials in $A$ such that at least one of them is of positive degree. Assume that there exists no large Zsigmondy prime for $(u, m)$, and that there exists a Zsigmondy prime for $(u, m)$. Then
\begin{itemize}

\item [(i)] $m + 1$ is the unique Zsigmondy prime for $(u, m)$;

\item [(ii)] either $\Psi_m(u) = \epsilon(m + 1)$ or $\Psi_m(u) = \epsilon \fq(m + 1)$ for some unit $\epsilon \in \bF_q^{\times}$ and some monic prime $\fq$ dividing $m$.

\end{itemize}

\end{corollary}

\begin{proof}

We contend that $\deg(m) > 0$. Indeed, assume the contrary, i.e., $\deg(m) = 0$, and hence $m = 1$. By assumption, $\deg(u) \ge 1$. Thus any monic prime $\wp$ dividing $u$ is a large Zsigmondy prime for $(u, m)$, which is a contradiction. Therefore $\deg(m) > 0$. 

Let $\fq$ be any Zsigmondy prime for $(u, m)$. By Corollary \ref{C-Luneburg-theorem-when-q->-2}$(i)$, $\fq$ does not divide $m$. Since $\fq$ is not a large Zsigmondy prime for $(u, m)$, we know from Remark \ref{R-large-Zsigmondy-primes} that $\deg(\fq) = \deg(m)$. Since $\fq$ is a Zsigmondy prime for $(u, m)$, we see from Remark \ref{Remark-The-remark-after-the-1st-definition-of-Zsigmondy-primes}$(i)$ that $m = P_{u, \fq}$. Hence it follows from Proposition \ref{Proposition-The-main-result-about-Carlitz-annihilator-of-u-and-wp}$(i)$ that $m$ divides $\fq - 1$, and thus $m = \fq - 1$ since $\deg(\fq - 1) = \deg(\fq) = \deg(m)$. Hence $\fq = m + 1$, which proves part $(i)$ of Corollary \ref{Corollary-Description-of-large-Zsigmondy-primes}.

By Remark \ref{Remark-A-Zsigmondy-prime-for-(u,m)-divides-Phi-m-(x)}$(i)$, we deduce from part $(i)$ that $m + 1$ divides $\Psi_m(u)$. We see that $(m + 1)^2$ does not divide $\Psi_m(u)$; otherwise, it follows from Proposition \ref{Proposition-The-2nd-basic-result-of-cyclotomic-polynomials} that $(m + 1)^2$ divides $C_m(u)$, and hence $m + 1$ is a large Zsigmondy prime for $(u, m)$, which is a contradiction. Therefore $\Psi_m(u)$ is of the form
\begin{align}
\label{Equation-The-equation-of-Phi-m-(u)-in-terms-of-Q-and-m+1-in-the-corollary-about-large-Zsigmondy-primes}
\Psi_m(u) = \epsilon Q(m + 1),
\end{align}
where $\epsilon \in \bF_q^{\times}$ and $Q$ is a monic polynomial such that $\gcd(Q, m + 1) = 1$. If $\deg(Q) = 0$, or equivalently $Q = 1$, we see that $\Psi_m(u) = \epsilon (m + 1)$, and part $(ii)$ follows. 

Suppose now that $\deg(Q) \ge 1$. Then there is a monic prime $\fq$ of positive degree dividing $Q$. By part $(i)$ and since $\gcd(\fq, m + 1) = 1$, we know that $\fq$ is a non-Zsigmondy prime for $(u, m)$, and it thus follows from Corollary \ref{C-Luneburg-theorem-when-q->-2}$(ii)$ that $m$ can be written in the form
\begin{align}
\label{Equation-The-equation-of-m-in-the-description-of-large-Zsigmondy-primes}
m = P_{u, \fq}\fq^s
\end{align}
for some positive integer $s$, where $P_{u, \fq}$ is the Carlitz annihilator of $(u, \fq)$. 

By Proposition \ref{Proposition-The-1st-basic-results-about-the-Carlitz-annihilator-of-(u, wp)}$(i)$, $P_{u, \fq}$ divides $\fq - 1$, which implies that $\deg(P_{u, \fq}) \le \deg(\fq - 1) = \deg(\fq)$. 

We now prove that $Q = \fq$. Assume the contrary, that is, $Q \ne \fq$. By Corollary \ref{C-Luneburg-theorem-when-q->-2}$(ii)$, $\fq^2$ does not divide $\Psi_m(u)$, and thus $Q/\fq$ is not divisible by $\fq$. Since $Q/\fq$ is a monic polynomial such that $Q/\fq \ne 1$, we deduce that $\deg(Q/\fq) > 0$, and thus there is a monic prime, say $\wp$, dividing $Q/\fq$. Since $\gcd(\fq, Q/\fq) = 1$, we deduce that $\gcd(\wp, \fq) = 1$. 

We will prove that $\wp = \fq - 1$ and $p = 2$. Indeed since $\wp$ divides $Q/\fq$, we deduce that $\gcd(\wp, m + 1) = 1$, and thus part $(i)$ tells us that $\wp$ is a non-Zsigmondy prime for $(u, m)$. Following the same arguments as above, one can write $m = P_{u, \wp}\wp^r$ for some $r \in \bZ_{>0}$, where $P_{u, \wp}$ is the Carlitz annihilator of $(u, \wp)$. Thus $m = P_{u, \fq}\fq^s = P_{u, \wp}\wp^r$.

By Proposition \ref{Proposition-The-1st-basic-results-about-the-Carlitz-annihilator-of-(u, wp)}$(i)$, $P_{u, \wp}$ divides $\wp - 1$, and thus $\deg(P_{u, \wp}) \le \deg(\wp - 1) = \deg(\wp)$. Since $\gcd(\wp, \fq) = 1$, we deduce that $\wp^r$ divides $P_{u, \fq}$, which in turn implies that $\wp^r$ divides $\fq - 1$, and hence $r\deg(\wp) \le \deg(\fq - 1) = \deg(\fq)$. Similarly, $\fq^s$ divides $\wp - 1$, and $s\deg(\fq) \le \deg(\wp - 1) = \deg(\wp)$. Therefore 
\begin{align*}
\deg(\wp) \le r\deg(\wp) \le \deg(\fq) \le s\deg(\fq) \le \deg(\wp),
\end{align*}
and thus
\begin{align*}
\deg(\wp) = r\deg(\wp) = \deg(\fq) = s\deg(\fq) = \deg(\wp).
\end{align*}
Hence $r = s = 1$, and  $\deg(\fq) = \deg(\wp)$. Since $\wp, \fq$ are monic, we deduce that 
\begin{align}
\label{Equation-The-equation-of-wp-in-the-description-of-large-Zsigmondy-primes}
\wp = \fq - 1,
\end{align}
and 
\begin{align}
\label{Equation-The-equation-of-fq-in-the-description-of-large-Zsigmondy-primes}
\fq = \wp - 1.
\end{align} 
Therefore $\fq = \wp - 1 = \fq - 2$, and hence $-2 = 0$, which implies that $p = 2$. (Recall that $p$ is the characteristic of $\bF_q$.)
 
Furthermore, since $\wp$ divides $P_{u, \fq}$, and $P_{u, \fq}$ divides $\fq - 1$, we deduce from (\ref{Equation-The-equation-of-wp-in-the-description-of-large-Zsigmondy-primes}) that 
\begin{align*}
P_{u, \fq} = \wp = \fq - 1,
\end{align*}
and thus 
\begin{align}
\label{Equation-The-equation-of-m-in-terms-of-q-in-the-corollary-about-large-Zsigmondy-primes}
m = P_{u, \fq}\fq = (\fq - 1)\fq.
\end{align}

In summary, we have shown that if $\wp$ is any monic prime dividing $Q/\fq$, then $\wp = \fq - 1$. In particular, this implies that $\fq - 1$ is a prime. Hence $Q/\fq = (\fq - 1)^{s_1}$ for some $s_1 \in \bZ_{>0}$. Since $\gcd(Q, m + 1) = 1$, we deduce from part $(i)$ that $\fq - 1$ is a non-Zsigmondy prime for $(u, m)$, and Corollary \ref{C-Luneburg-theorem-when-q->-2}$(ii)$ tells us that $s_1 = 1$. Hence $Q = \fq(\fq - 1)$, and It therefore follows from (\ref{Equation-The-equation-of-Phi-m-(u)-in-terms-of-Q-and-m+1-in-the-corollary-about-large-Zsigmondy-primes}) and (\ref{Equation-The-equation-of-m-in-terms-of-q-in-the-corollary-about-large-Zsigmondy-primes}) that
\begin{align}
\label{Equation-The-equation-of-Psi-m-(u)-in-terms-of-m-only-in-the-corollary-about-large-Zsigmondy-primes}
\Psi_m(u) = \epsilon Q(m + 1) = \epsilon \fq(\fq - 1)(m + 1) = \epsilon m(m + 1).
\end{align}

We consider two cases:

$\star$ \textit{Case 1. $\deg(u) = 0$.}

We see that $u = 1$. From (\ref{Equation-The-equation-of-m-in-terms-of-q-in-the-corollary-about-large-Zsigmondy-primes}) and (\ref{Equation-The-equation-of-Psi-m-(u)-in-terms-of-m-only-in-the-corollary-about-large-Zsigmondy-primes}), we have $\Psi_m(u) = \epsilon (\fq - 1)\fq(\fq(\fq - 1) + 1)$, and thus 
\begin{align}
\label{Equation-The-degree-of-Psi-m-(u)-in-terms-of-fq-in-the-corollary-about-large-Zsigmondy-primes}
\deg(\Psi_m(u)) = 4\deg(\fq).
\end{align}

On the other hand, since $\deg(u) = 0$, Lemma \ref{Lemma-The-degrees-of-cyclotomic-polynomials-when-u-is-a-unit} yields
\begin{align}
\label{Equation-The-degree-of-Psi-m-(u)-in-terms-of-Phi(m)-large-Zsigmondy-primes-corollary}
\deg(\Psi_m(u)) = \dfrac{\Phi(m) + \delta}{q}
\end{align}
for some integer $\delta \in \{-1, 0, 1\}$. Applying Lemma \ref{Lemma-The-lower-bound-of-Phi(m)}, we see from (\ref{Equation-The-equation-of-m-in-terms-of-q-in-the-corollary-about-large-Zsigmondy-primes}) that
\begin{align*}
\Phi(m) = \Phi(\fq(\fq - 1)) = \Phi(\fq)\Phi(\fq - 1) \ge (q - 1)^2\deg(\fq)\deg(\fq - 1) = (q - 1)^2\deg(\fq)^2.
\end{align*}
Since $\delta \ge -1$, we deduce from (\ref{Equation-The-degree-of-Psi-m-(u)-in-terms-of-fq-in-the-corollary-about-large-Zsigmondy-primes}) and (\ref{Equation-The-degree-of-Psi-m-(u)-in-terms-of-Phi(m)-large-Zsigmondy-primes-corollary}) that
\begin{align*}
4\deg(\fq) = \deg(\Psi_m(u)) \ge \dfrac{\Phi(m) - 1}{q} \ge \dfrac{(q - 1)^2\deg(\fq)^2 - 1}{q},
\end{align*}
and thus
\begin{align}
\label{Inequality-The-1st-inequality-large-Zsigmondy-primes-corollary}
\deg(\fq)((q - 1)^2\deg(\fq) - 4q) \le 1.
\end{align}

Note that $q \ge 4$ since $q$ is a power of 2 and $q > 2$. If $\deg(\fq) \ge 2$, then we see that 
\begin{align*}
2(q - 1)^2 - 4q = 2(q^2 - 4q + 1) = 2(q(q - 4) + 1) \ge 2,
\end{align*}
and thus $\deg(\fq)((q - 1)^2\deg(\fq) - 4q) \ge 4$, which is a contradiction to (\ref{Inequality-The-1st-inequality-large-Zsigmondy-primes-corollary}). 

Suppose now that $\deg(\fq) = 1$. Since $q$ is a power of $2$ and $q > 2$, either $q \ge 8$ or $q = 4$. If $q = 4$, then we see from (\ref{Equation-The-equation-of-m-in-terms-of-q-in-the-corollary-about-large-Zsigmondy-primes}) that $m = \fq(\fq - 1)$, where $\fq$ is a monic prime of degree 1 in $\bF_{2^2}[T]$. Since $u = 1$, Theorem \ref{Theorem-Zsigmondy-theorem-in-characteristic-two} tells us that there exist no Zsigmondy primes for $(u, m)$ in $\bF_{2^2}[T]$, which is a contradiction. Thus $q \ge 8$, and therefore
\begin{align*}
\deg(\fq)((q - 1)^2\deg(\fq) - 4q) = q^2 - 6q + 1 = q(q - 6) + 1 \ge 17,
\end{align*}
which is a contradiction to (\ref{Inequality-The-1st-inequality-large-Zsigmondy-primes-corollary}).

$\star$ \textit{Case 2. $\deg(u) > 0$.}

By Lemma \ref{Lemma-The-degrees-of-cyclotomic-polynomials}, (\ref{Equation-The-equation-of-m-in-terms-of-q-in-the-corollary-about-large-Zsigmondy-primes}), and (\ref{Equation-The-equation-of-Psi-m-(u)-in-terms-of-m-only-in-the-corollary-about-large-Zsigmondy-primes}), we know that $\deg(\Psi_m(u)) = \Phi(m)\deg(u) = 2\deg(m)$. Thus, by Lemma \ref{Lemma-The-lower-bound-of-Phi(m)}, $2\deg(m) = \Phi(m)\deg(u) \ge \Phi(m) \ge (q - 1)\deg(m)$, and therefore $(q - 3)\deg(m) \le 0$, which is a contradiction since $q \ge 4$ and $\deg(m) \ge 1$. 

By \textit{Cases 1 and 2}, we conclude that $Q = \fq$, and it thus follows from (\ref{Equation-The-equation-of-Phi-m-(u)-in-terms-of-Q-and-m+1-in-the-corollary-about-large-Zsigmondy-primes}) that $\Psi_m(u) =  \epsilon \fq (m + 1)$, which proves part $(ii)$.

\end{proof}

We now state several lemmas (see Lemmas \ref{Lemma-The-exceptional-case-EC-III}, \ref{Lemma-The-exceptional-case-EC-IV}, \ref{Lemma-The-exceptional-case-EC-V}, \ref{Lemma-The-exceptional-case-VI}, \ref{L-the-exceptional-case-EC-VI}, \ref{L-the-exceptional-case-EC-VII}, \ref{L-the-exceptional-case-EC-VIII}, and \ref{L-the-exceptional-case-EC-IX}) that we need in the proof of Theorem \ref{Theorem-The-main-theorem-about-large-Zsigmondy-primes}. These lemmas rule out the exceptional cases in Theorem \ref{Theorem-The-main-theorem-about-large-Zsigmondy-primes} that naturally appear in the proof of Theorem \ref{Theorem-The-main-theorem-about-large-Zsigmondy-primes}. The proofs of the lemmas are purely computational, and can be easily verified with the aid of a computer algebra system. So we do not include the proofs here.

\begin{lemma}
\label{Lemma-The-exceptional-case-EC-III}

Let $q = 3$. Let $\cX_3$ be the set of all monic polynomials $m \in \bF_3[T]$ satisfying the following conditions:
\begin{itemize}

\item [(i)] $m = (\wp - 1)\wp^s$ for some monic prime $\wp$ of degree one in $\bF_3[T]$ and some integer $2 \le s \le 7$; 

\item [(ii)] $m + 1$ is the only Zsigmondy prime for $(1, m)$;

\item [(iii)] there are no large Zsigmondy primes for $(1, m)$. 

\end{itemize}
Then $\cX_3 =\{(T - 1)T^2, T(T + 1)^2, (T + 1)(T + 2)^2\}.$

\end{lemma}

\begin{lemma}
\label{Lemma-The-exceptional-case-EC-IV}

Let $q = 3$. Let $\cX_4$ be the set of all polynomials $m \in \bF_3[T]$ satisfying the following conditions:
\begin{itemize}

\item [(i)] $m$ is a monic prime of degree one in $\bF_3[T]$; 

\item [(ii)] $m + 1$ is the only Zsigmondy prime for $(m, m)$;

\item [(iii)] there are no large Zsigmondy primes for $(m, m)$. 

\end{itemize}
Then $\cX_4 = \{T, T + 1, T + 2\}.$

\end{lemma}

\begin{lemma}
\label{Lemma-The-exceptional-case-EC-V}

Let $q = 3$. Let $\cX_5$ be the set of all monic polynomials $m \in \bF_3[T]$ satisfying the following two conditions:
\begin{itemize}

\item [(i)] there exists a monic prime $\fq$ of degree 2 or 3 in $\bF_3[T]$ such that the Carlitz annihilator $P_{1, \fq}$ of $(1, \fq)$ is of degree 2 and $m = P_{1, \fq}\fq$;

\item [(ii)] $m + 1$ is a prime in $\bF_3[T]$.

\end{itemize}
Then $\cX_5 = \emptyset$.

\end{lemma}

\begin{lemma}
\label{Lemma-The-exceptional-case-VI}

Let $q = 5$. Let $\cX_6$ be the set of all monic polynomials $m \in \bF_5[T]$ satisfying the following conditions: 
\begin{itemize}

\item [(i)] $m = (\wp - 1)\wp$ for some monic prime $\wp$ of degree one in $\bF_5[T]$; 

\item [(ii)] $m + 1$ is the only Zsigmondy prime for $(1, m)$; and

\item [(iii)] there are no large Zsigmondy primes for $(1, m)$.

\end{itemize}

Then
\begin{align*}
\cX_6 &= \{m = (\wp - 1)\wp \; | \; \text{$\wp$ is a monic prime of degree one in $\bF_5[T]$}\} \\
&= \{m = (T + \alpha - 1)(T + \alpha) \; |\; \alpha \in \bF_5\}.
\end{align*}

\end{lemma}

\begin{lemma}
\label{L-the-exceptional-case-EC-VI}

Assume that $q > 2$. Let $\cX_7$ be the set of all monic primes of degree one in $\bF_q[T]$. Then there are no large Zsigmondy primes for $(1, m)$ for any $m \in \cX_7$.

\end{lemma}

\begin{lemma}
\label{L-the-exceptional-case-EC-VII}

Let $q = 3$. Let $\cX_8$ be the set of all monic polynomials $m \in \bF_3[T]$ satisfying the following conditions:
\begin{itemize}

\item [(i)] $m = \wp^2$ for some monic prime $\wp$ of degree one in $\bF_3[T]$;

\item [(ii)] $m + 1$ is the only Zsigmondy prime for $(1, m)$;

\item [(iii)] there are no large Zsigmondy primes for $(1, m)$.

\end{itemize}
Then $\cX_8 = \{T^2, (T + 1)^2, (T + 2)^2\}.$

\end{lemma}

\begin{lemma}
\label{L-the-exceptional-case-EC-VIII}

Let $q = 4$, and write $\bF_4 = \bF_2(w)$, where $w$ is an element in the algebraic closure of $\bF_2$ such that $w^2 + w + 1 = 0$. Let $\cX_9$ be the set of all monic polynomials $m \in \bF_4[T]$ satisfying the following conditions:
\begin{itemize}

\item [(i)] $m = m_1m_2$, where $m_1, m_2$ are monic polynomials in $\bF_4[T]$ such that $\deg(m_1) = \deg(m_2) = 1$ and $\gcd(m_1, m_2) = 1$;

\item [(ii)] $m + 1$ is the only Zsigmondy prime for $(1, m)$;

\item [(iii)] there are no large Zsigmondy primes for $(1, m)$.

\end{itemize}
Then $\cX_9 = \{T(T + w), T(T + w^2), (T + 1)(T + w), (T + 1)(T + w^2)\}.$

\end{lemma}

\begin{lemma}
\label{L-the-exceptional-case-EC-IX}

Let $q = 3$. Let $\cX_{10}$ be the set of all monic polynomials $m \in \bF_3[T]$ satisfying the following conditions:
\begin{itemize}

\item [(i)] $m$ is square-free, i.e., $\wp^2$ does not divide $m$ for any monic prime $\wp$;

\item [(ii)] $m = m_1m_2$, where $m_1, m_2$ are monic polynomials in $\bF_3[T]$ such that $\gcd(m_1, m_2) = 1$, $\deg(m_1) \in \{1, 2\}$, and $\deg(m_2) = 1$;

\item [(iii)] $m + 1$ is the only Zsigmondy prime for $(1, m)$;

\item [(iv)] there are no large Zsigmondy primes for $(1, m)$.

\end{itemize}
Then $\cX_{10} = \{T^3 + 2T\}.$

\end{lemma}

In order to rule out some exceptional cases in the proof of Theorem \ref{Theorem-The-main-theorem-about-large-Zsigmondy-primes}, we need to strengthen Lemma \ref{Lemma-The-lower-bound-of-Phi(m)}, and obtain a sharper lower bound for $\Phi(m)$ for some special cases of $m$.

\begin{lemma}
\label{L-lower-bound-for-Phi(prime-of-deg->=2)}

Assume that $q > 2$. Let $\wp$ be a monic prime of degree $\ge 2$. Then
\begin{align}
\label{E-lower-bound-for-Phi(prime-of-deg->=2}
\Phi(\wp) - q\deg(\wp) \ge q^2 - 2q - 1.
\end{align}
Furthermore equality in (\ref{E-lower-bound-for-Phi(prime-of-deg->=2}) occurs if and only if $\deg(\wp) = 2$. 

\end{lemma}

\begin{proof}

Let $H(\alpha)$ be the function defined by
\begin{align}
\label{E-the-function-H(alpha)-in-L-about-lower-bound-for-Phi(prime-of-deg->=2}
H(\alpha) = q^{\alpha} - 1 - q\alpha,
\end{align}
where $\alpha$ ranges over the set $[2, \infty)$. Note that (\ref{E-lower-bound-for-Phi(prime-of-deg->=2}) is equivalent to the inequality
\begin{align*}
H(\deg(\wp)) \ge q^2 - 2q - 1.
\end{align*}

We prove that $H$ is a strictly increasing function over the interval $[2, \infty)$. Indeed, the derivative of $H$ is equal to
\begin{align*}
H^{\prime}(\alpha) = \ln(q)q^{\alpha} - q,
\end{align*}
and since $q \ge 3$ and $\alpha \ge 2$, we deduce that
\begin{align*}
H^{\prime}(\alpha) \ge \ln(3)q^2 - q > q^2 - q = q(q - 1) \ge 6 > 0.
\end{align*}
Thus $H$ is a strictly increasing function over the interval $[2, \infty)$, and therefore
\begin{align}
\label{E-inequality-of-H(alpha)-in-L-about-lower-bound-for-Phi(prime-of-deg->=2}
H(\alpha) \ge H(2) = q^2 - 2q - 1.
\end{align} 
Note that equality in (\ref{E-inequality-of-H(alpha)-in-L-about-lower-bound-for-Phi(prime-of-deg->=2}) occurs if and only if $\alpha = 2$. 

Now replacing $\alpha$ by $\deg(\wp)$ in (\ref{E-inequality-of-H(alpha)-in-L-about-lower-bound-for-Phi(prime-of-deg->=2}), Lemma \ref{L-lower-bound-for-Phi(prime-of-deg->=2)} follows immediately.

\end{proof}

Since $q^2 - 2q - 1 = q(q - 2) - 1 \ge 2$, the next result follows immediately from Lemma \ref{L-lower-bound-for-Phi(prime-of-deg->=2)}.

\begin{corollary}
\label{C-lower-bound-for-Phi(prime-of-deg->=2)}

Assume that $q > 2$. Let $\wp$ be a monic prime of degree $\ge 2$. Then $\Phi(\wp) > q\deg(\wp).$

\end{corollary}

\begin{lemma}
\label{L-lower-bound-for-Phi(prime-power)}

Assume that $q > 2$. Let $\wp$ be a monic prime, and let $s$ be an integer such that $s \ge 2$. Then
\begin{align}
\label{E-lower-bound-for-Phi(prime-power)}
\Phi(\wp^s) - q\deg(\wp^s) \ge q(q - 3)
\end{align}
Furthermore equality in (\ref{E-lower-bound-for-Phi(prime-power)}) occurs if and only if $\deg(\wp) = 1$ and $s = 2$.

\end{lemma}

\begin{proof}

Let $F(\alpha)$ be the function defined by $F(\alpha) = q^{s\alpha} - q^{(s - 1)\alpha} - sq\alpha$, where $\alpha$ ranges over the set $[1, \infty)$. Note that (\ref{E-lower-bound-for-Phi(prime-power)}) is equivalent to the inequality $F(\deg(\wp)) \ge q(q - 3)$. In order to prove Lemma \ref{L-lower-bound-for-Phi(prime-power)}, it suffices to show that $F$ is a strictly increasing function over the interval $[1, \infty)$. For the proof of the latter, one can use the same arguments as in the proof of Lemma \ref{L-lower-bound-for-Phi(prime-of-deg->=2)}.

\end{proof}

The next result is immediate from Lemma \ref{L-lower-bound-for-Phi(prime-power)}.

\begin{corollary}
\label{C-the-weaker-bound-for-Phi(wp^s)}

Assume that $q > 2$. Let $\wp$ be a monic prime, and let $s$ be an integer such that $s \ge 2$. Then
\begin{align}
\label{E-the-weaker-lower-bound-for-Phi(wp^s)}
\Phi(\wp^s) \ge q\deg(\wp^s).
\end{align}
Furthermore equality in (\ref{E-the-weaker-lower-bound-for-Phi(wp^s)}) occurs if and only if $q = 3$, $\deg(\wp) = 1$, and $s = 2$. 

\end{corollary}

The following theorem is the second main result of this paper that can be viewed as a function field analogue of Feit's theorem (see Feit \cite[Theorem A]{Feit-PAMS-1988}). Corollary \ref{Corollary-Description-of-large-Zsigmondy-primes} plays a key role in the proof of the next theorem.

\begin{theorem}
\label{Theorem-The-main-theorem-about-large-Zsigmondy-primes}

Assume that $q > 2$. Let $m, u$ be monic polynomials in $A$ such that at least one of them is of positive degree. Then there exists a large Zsigmondy prime for $(u, m)$ except exactly in the following cases:
\begin{itemize}

\item [(EC-I)] $q = 3, u = 1$, and $m = (\wp - 1)\wp$, where $\wp$ is an arbitrary monic prime of degree one in $\bF_3[T]$;

\item [(EC-II)] $q = 2^2, u = 1$, and $m = (\wp - 1)\wp$, where $\wp$ is an arbitrary monic prime of degree one in $\bF_{2^2}[T]$;

\item [(EC-III)] $q = 3, u = 1$, and 
\begin{align*}
m \in \cX_3 = \{(T - 1)T^2, T(T + 1)^2, (T + 1)(T + 2)^2\},
\end{align*}
where $\cX_3$ is the set in Lemma \ref{Lemma-The-exceptional-case-EC-III}.

\item [(EC-IV)] $q = 3$, and
\begin{align*}
u = m \in \cX_4 = \{T, T + 1, T + 2\},
\end{align*}  
where $\cX_4$ is the set in Lemma \ref{Lemma-The-exceptional-case-EC-IV}.

\item [(EC-V)] $q = 5, u = 1$, and
\begin{align*}
m \in \cX_6 = \{T(T + 1), (T + 1)(T + 2), (T + 2)(T + 3), (T + 3)(T + 4), (T + 4)T\},
\end{align*}
where $\cX_6$ is the set in Lemma \ref{Lemma-The-exceptional-case-VI}.

\item [(EC-VI)] $u = 1$, and $m \in \cX_7$, i.e., $m$ is a monic prime of degree one in $\bF_q[T]$, where $\cX_7$ is the set in Lemma \ref{L-the-exceptional-case-EC-VI}. (Note that there is no restriction on $q$ in this exceptional case.)

\item [(EC-VII)] $q = 3$, $u = 1$, and 
\begin{align*}
m \in \cX_8 = \{T^2, (T + 1)^2, (T + 2)^2\}, 
\end{align*}
where $\cX_8$ is the set in Lemma \ref{L-the-exceptional-case-EC-VII}.

\item [(EC-VIII)] $q = 4$, $u = 1$, and 
\begin{align*}
m \in \cX_9 = \{T(T + w), T(T + w^2), (T + 1)(T + w), (T + 1)(T + w^2)\}, 
\end{align*}
where $\cX_9$ is the set in Lemma \ref{L-the-exceptional-case-EC-VIII}.(Note that $\bF_4 = \bF_2(w)$, where $w$ is an element in the algebraic closure of $\bF_2$ such that $w^2 + w + 1 = 0$.)

\item [(EC-IX)] $q = 3$, $u = 1$, and 
\begin{align*}
m \in \cX_{10} = \{T^3 + 2T\}, 
\end{align*}
where $\cX_{10}$ is the set in Lemma \ref{L-the-exceptional-case-EC-IX}.

\end{itemize}

\end{theorem}

\begin{proof}

If $\deg(m) = 0$, then $m = 1$, and it thus follows from the assumption that $u$ is of positive degree. It is easy to see that any monic prime $\wp$ dividing $u$ is a large Zsigmondy prime for $(u, m) = (u, 1)$. For the rest of the proof, without loss of generality, one can assume that $m$ is of positive degree.

We first consider the cases when we are in one of the exceptional cases (EC-I)--(EC-IX). If we are in the exceptional case (EC-I) or the exceptional case (EC-II), then Theorem \ref{Theorem-The-full-general-version-of-Zsigmondy-Theorem-with-p-not--equal-to-2} and Theorem \ref{Theorem-Zsigmondy-theorem-in-characteristic-two} tell us that there are no Zsigmondy primes for $(u, m)$, and thus there are no large Zsigmondy primes for $(u, m)$.

On the other hand, Lemmas \ref{Lemma-The-exceptional-case-EC-III}, \ref{Lemma-The-exceptional-case-EC-IV}, \ref{Lemma-The-exceptional-case-VI}, \ref{L-the-exceptional-case-EC-VI}, \ref{L-the-exceptional-case-EC-VII}, \ref{L-the-exceptional-case-EC-VIII}, and \ref{L-the-exceptional-case-EC-IX} tell us that there are no large Zsigmondy primes for $(u, m)$ if we are in one of the exceptional cases (EC-III)--(EC-IX). 

Suppose, for the rest of the proof, that we are not in any of the exceptional cases (EC-I)--(EC-IX). We prove that there exists a large Zsigmondy prime for $(u, m)$. Assume the contrary, that is, 
\begin{itemize}

\item [(LZP0)] \textit{there exist no large Zsigmondy primes for $(u, m)$.}

\end{itemize}
By Theorem \ref{Theorem-The-full-general-version-of-Zsigmondy-Theorem-with-p-not--equal-to-2} and Theorem \ref{Theorem-Zsigmondy-theorem-in-characteristic-two}, we know that there exists a Zsigmondy prime for $(u, m)$, and it thus follows from Corollary \ref{Corollary-Description-of-large-Zsigmondy-primes} that the following are true:
\begin{itemize}

\item [(LZP1)] \textit{$m + 1$ is the only Zsigmondy prime for $(u, m)$; and}

\item [(LZP2)] \textit{either $\Psi_m(u) = \epsilon (m + 1)$ for some unit $\epsilon \in \bF_q^{\times}$ or $\Psi_m(u) = \epsilon \fq(m + 1)$ for some unit $\epsilon \in \bF_q^{\times}$ and some monic prime $\fq$ dividing $m$.}

\end{itemize}

We consider the following two cases:

$\star$ \textit{Case 1. $\Psi_m(u) = \epsilon \fq(m + 1)$ for some unit $\epsilon \in \bF_q^{\times}$ and some monic prime $\fq$ dividing $m$.}

Since $\gcd(m, m + 1) = 1$ and $\fq$ divides $m$, we deduce that $\gcd(\fq, m + 1) = 1$. Hence we, by appealing to (LZP1), find that $\fq$ is a non-Zsigmondy prime for $(u, m)$. Since $\fq$ divides $\Psi_m(u)$, it follows from Corollary \ref{C-Luneburg-theorem-when-q->-2}$(ii)$ that $m$ is of the form
\begin{align}
\label{Equation-Eqn-of-m=P_{u-fq}-fq^s--in-Case-1-LZP-the-main-theorem}
m = P_{u, \fq}\fq^s,
\end{align}
where $P_{u, \fq}$ is the Carlitz annihilator of $(u, \fq)$ and $s$ is a positive integer. 

We consider the following two subcases, according as to whether $\deg(u) \ge 1$ or $\deg(u) = 0$. 

$\star$ \textit{Subcase 1A. $\deg(u) \ge 1$.} 

Since $\deg(P_{u, \fq}) \ge 0$, we consider the following two subsubcases, according as to whether $\deg(P_{u, \fq}) = 0$ or $\deg(P_{u, \fq}) \ge 1$.

$\bullet$ \textit{Subsubcase 1A(i). $\deg(P_{u, \fq}) = 0$.}

By Definition \ref{Definition-The-Carlitz-annihilator-of-(u, wp)}, we see that $P_{u, \fq} = 1$. From (\ref{Equation-Eqn-of-m=P_{u-fq}-fq^s--in-Case-1-LZP-the-main-theorem}), we find that
\begin{align}
\label{E-eqn-of-m=fq^s--in-Subcase-1A-of-Case-1-in-LZP-theorem}
m = \fq^s.
\end{align}

Since $\deg(u) \ge 1$, we deduce from (\ref{E-eqn-of-m=fq^s--in-Subcase-1A-of-Case-1-in-LZP-theorem}), Lemma \ref{Lemma-The-degrees-of-cyclotomic-polynomials}, and Lemma \ref{Lemma-The-lower-bound-of-Phi(m)} that 
\begin{align}
\label{E-inequality-of-deg(Psi_m(u))-ge-s(q-1)deg(fq)--in--Subsubcase-1A(i)-of-Subcase-1A-in-Case-1--LZP-thm}
\deg(\Psi_m(u)) = \Phi(m)\deg(u) \ge (q - 1)\deg(m) = s(q - 1)\deg(\fq).
\end{align}
Note that equality in (\ref{E-inequality-of-deg(Psi_m(u))-ge-s(q-1)deg(fq)--in--Subsubcase-1A(i)-of-Subcase-1A-in-Case-1--LZP-thm}) occurs if and only if $\deg(u) = 1$ and $\Phi(m) = (q - 1)\deg(m)$.

On the other hand, we see from (\ref{E-eqn-of-m=fq^s--in-Subcase-1A-of-Case-1-in-LZP-theorem}) that
\begin{align}
\label{E-deg(Psi_m(u))=(s+1)deg(fq)---in--Subsubcase-1A(i)--of--Subcase-1A--in-Case-1-LZP-thm}
\deg(\Psi_m(u)) = \deg(\epsilon \fq(m + 1)) = \deg(\fq) + \deg(m + 1) = \deg(\fq) + \deg(m) = (s + 1)\deg(\fq).
\end{align}
Combining (\ref{E-inequality-of-deg(Psi_m(u))-ge-s(q-1)deg(fq)--in--Subsubcase-1A(i)-of-Subcase-1A-in-Case-1--LZP-thm}), (\ref{E-deg(Psi_m(u))=(s+1)deg(fq)---in--Subsubcase-1A(i)--of--Subcase-1A--in-Case-1-LZP-thm}), we find that $(s + 1)\deg(\fq) = \deg(\Psi_m(u)) \ge s(q - 1)\deg(\fq)$, and thus
\begin{align}
\label{E--(s(q-1)-(s+1))deg(fq)--le--0---the-1st-key-inequality---in--Subsubcase--1A(i)--of--Subcase-1A--in--Case-1-LZP-thm}
(s(q - 1) - (s + 1))\deg(\fq) \le 0.
\end{align}
By the remark following (\ref{E-inequality-of-deg(Psi_m(u))-ge-s(q-1)deg(fq)--in--Subsubcase-1A(i)-of-Subcase-1A-in-Case-1--LZP-thm}), note that equality in (\ref{E--(s(q-1)-(s+1))deg(fq)--le--0---the-1st-key-inequality---in--Subsubcase--1A(i)--of--Subcase-1A--in--Case-1-LZP-thm}) occurs if and only if $\deg(u) = 1$ and $\Phi(m) = (q - 1)\deg(m)$.

Furthermore since $s \ge 1$ and $q \ge 3$, we find that $s(q - 1) - (s + 1) \ge 2s - (s + 1) = s - 1 \ge 0$, and since $\deg(\fq) \ge 1$, we deduce that
\begin{align}
\label{E--inequality--(s(q-1)-(s+1))deg(fq)--ge--0--the-2nd-key-inequality-in-Subsubcase-1A(i)--of--Subcase-1A-in-Case-1-LZP-thm}
(s(q - 1) - (s + 1))\deg(\fq) \ge 0.
\end{align}
Note that equality in (\ref{E--inequality--(s(q-1)-(s+1))deg(fq)--ge--0--the-2nd-key-inequality-in-Subsubcase-1A(i)--of--Subcase-1A-in-Case-1-LZP-thm}) occurs if and only if $s = 1$ and $q = 3$. 

Combining (\ref{E--(s(q-1)-(s+1))deg(fq)--le--0---the-1st-key-inequality---in--Subsubcase--1A(i)--of--Subcase-1A--in--Case-1-LZP-thm}), (\ref{E--inequality--(s(q-1)-(s+1))deg(fq)--ge--0--the-2nd-key-inequality-in-Subsubcase-1A(i)--of--Subcase-1A-in-Case-1-LZP-thm}), we deduce that $(s(q - 1) - (s + 1))\deg(\fq) = 0$. Hence by appealing to the remarks following (\ref{E--(s(q-1)-(s+1))deg(fq)--le--0---the-1st-key-inequality---in--Subsubcase--1A(i)--of--Subcase-1A--in--Case-1-LZP-thm}), (\ref{E--inequality--(s(q-1)-(s+1))deg(fq)--ge--0--the-2nd-key-inequality-in-Subsubcase-1A(i)--of--Subcase-1A-in-Case-1-LZP-thm}), we deduce that $s = 1$, $q = 3$, $\deg(u) = 1$ and $\Phi(m) = (q - 1)\deg(m)$. It then follows from (\ref{E-eqn-of-m=fq^s--in-Subcase-1A-of-Case-1-in-LZP-theorem}) that $m = \fq$, and 
\begin{align}
\label{E-eqn-Phi(fq)=2deg(fq)--in-Subsubcase-1A(i)--of--Subcase--1A--in--Case-1-LZP-thm}
\Phi(\fq) = \Phi(m) = 2\deg(m) = 2\deg(\fq).
\end{align}

If $\deg(\fq) \ge 2$, then we deduce from Corollary \ref{C-lower-bound-for-Phi(prime-of-deg->=2)} that $\Phi(\fq) > q\deg(\fq) = 3\deg(\fq)$, which is a contradiction to (\ref{E-eqn-Phi(fq)=2deg(fq)--in-Subsubcase-1A(i)--of--Subcase--1A--in--Case-1-LZP-thm}). Hence we deduce that $\deg(\fq) = 1$. 

Recall that $\deg(u) = 1$ and $P_{u, \fq} = 1$. By Proposition \ref{Proposition-The-1st-basic-results-about-the-Carlitz-annihilator-of-(u, wp)}$(ii)$, we find that $\fq$ divides $u$, and since $\deg(u) = \deg(\fq) = 1$ and $\fq, u$ are monic polynomials, we deduce that $\fq = u$. 

In summary, we have showed that $q = 3$, and $m = u = \fq$ which is a monic prime of degree one in $\bF_3[T]$. By appealing to (LZP0) and (LZP1), this implies that we are in the exceptional case (EC-IV), which is a contradiction.

$\bullet$ \textit{Subsubcase 1A(ii). $\deg(P_{u, \fq}) \ge 1$.}

Set
\begin{align}
\label{E-lambda=deg(P_{u--fq})--in-Subsubcase-1A(ii)---of-Subcase-1A--in--Case-1--LZP-thm}
\lambda = \deg(P_{u, \fq}) \ge 1,
\end{align}
and 
\begin{align}
\label{E-gamma=deg(fq)--in-Subsubcase-1A(ii)---of-Subcase-1A--in--Case-1--LZP-thm}
\gamma = \deg(\fq) \ge 1.
\end{align}

Recall from Proposition \ref{Proposition-The-1st-basic-results-about-the-Carlitz-annihilator-of-(u, wp)}$(i)$ that $P_{u, \fq}$ divides $\fq - 1$, and since $\gcd(\fq, \fq - 1) = 1$, we deduce that $\gcd(P_{u, \fq}, \fq) = 1$. Hence by appealing to Lemma \ref{Lemma-The-lower-bound-of-Phi(m)}, we deduce from (\ref{Equation-Eqn-of-m=P_{u-fq}-fq^s--in-Case-1-LZP-the-main-theorem}) that
\begin{align}
\label{E-inequality-Phi(m)--ge--(q-1)^2deg(P_{u--fq})deg(fq)---in--Subsubcase-1A(ii)--of---Subcase-1A-in---Case--1---LZP--thm}
\Phi(m) = \Phi(P_{u, \fq} \fq^s) = \Phi(P_{u, \fq})\Phi(\fq^s) \ge (q - 1)^2\deg(P_{u, \fq})\deg(\fq^s) = s(q - 1)^2\deg(P_{u, \fq})\deg(\fq).
\end{align}
Since $\deg(u) \ge 1$, we deduce from Lemma \ref{Lemma-The-degrees-of-cyclotomic-polynomials} and (\ref{E-inequality-Phi(m)--ge--(q-1)^2deg(P_{u--fq})deg(fq)---in--Subsubcase-1A(ii)--of---Subcase-1A-in---Case--1---LZP--thm}) that
\begin{align}
\label{E--inequality-of-deg(Psi_m(u))---ge---(q-1)^2deg(P_{u--fq})deg(fq)---in--Subsubcase-1A(ii)--of---Subcase-1A--in-Case-1-LZP-thm}
\deg(\Psi_m(u)) = \deg(u)\Phi(m) \ge \Phi(m) \ge s(q - 1)^2\deg(P_{u, \fq})\deg(\fq) = s(q - 1)^2\lambda\gamma.
\end{align}
 
 Recall that $\Psi_m(u) = \epsilon \fq(m + 1)$, and it thus follows from (\ref{Equation-Eqn-of-m=P_{u-fq}-fq^s--in-Case-1-LZP-the-main-theorem}) that
 \begin{align}
 \label{E-eqn-of-deg(Psi_m(u))--==---(s+1)gamma+lambda--in-Subsubcase-1A(ii)--in-Subcase-1A-in-Case-1---LZP-thm}
 \deg(\Psi_m(u)) = \deg(\fq) + \deg(m + 1) &= \deg(\fq) + \deg(m) \nonumber \\
 &= (s + 1)\deg(\fq) + \deg(P_{u, \fq}) \nonumber \\
 &= (s + 1)\gamma + \lambda.
 \end{align} 
 
 Combining (\ref{E--inequality-of-deg(Psi_m(u))---ge---(q-1)^2deg(P_{u--fq})deg(fq)---in--Subsubcase-1A(ii)--of---Subcase-1A--in-Case-1-LZP-thm}) and (\ref{E-eqn-of-deg(Psi_m(u))--==---(s+1)gamma+lambda--in-Subsubcase-1A(ii)--in-Subcase-1A-in-Case-1---LZP-thm}), we find that
 \begin{align*}
 (s + 1)\gamma + \lambda = \deg(\Psi_m(u)) \ge s(q - 1)^2\lambda\gamma,
 \end{align*}
 and thus
 \begin{align}
 \label{E-eqn-of-s(q-1)^2lambda-gamma--(s+1)gamma--lambda--le-0-Subsubcase-1A(ii)-Subcase-1A-Case-1-LZP-thm}
 s(q - 1)^2\lambda\gamma -  (s + 1)\gamma - \lambda \le 0.
 \end{align}
  
 Since $\lambda \ge 1$, $\gamma \ge 1$, and $s \ge 1$, we deduce that
 \begin{align}
 \label{E-eqn-of-2slambda--gamma--ge---(s-1)gamma--ge--0-Subsubcase-1A(ii)-Subcase-1A-Case-1-LZP-thm}
 2s\lambda\gamma - (s + 1)\gamma \ge 2s\gamma - (s + 1)\gamma = (s - 1)\gamma \ge 0.
 \end{align}
 On the other hand, since $\lambda \ge 1$, $\gamma \ge 1$, and $s \ge 1$, we see that
 \begin{align}
 \label{E-eqn-2s-lambda--gamma--ge---1--in--Subsubcase-1A(ii)-Subcase-1A-Case-1-LZP-thm}
 2s\lambda\gamma - \lambda \ge 2\lambda - \lambda = \lambda \ge 1.
 \end{align}
 Combining (\ref{E-eqn-of-2slambda--gamma--ge---(s-1)gamma--ge--0-Subsubcase-1A(ii)-Subcase-1A-Case-1-LZP-thm}) and (\ref{E-eqn-2s-lambda--gamma--ge---1--in--Subsubcase-1A(ii)-Subcase-1A-Case-1-LZP-thm}), and note that $(q - 1)^2 \ge 4$, we deduce that
 \begin{align*}
 s(q - 1)^2\lambda\gamma -  (s + 1)\gamma - \lambda &\ge 4s\lambda\gamma -  (s + 1)\gamma - \lambda \\
 &= (2s\lambda\gamma -  (s + 1)\gamma) + (2s\lambda\gamma -\lambda )\\
 &\ge 0 + 1 = 1,
 \end{align*}
 which is a contradiction to (\ref{E-eqn-of-s(q-1)^2lambda-gamma--(s+1)gamma--lambda--le-0-Subsubcase-1A(ii)-Subcase-1A-Case-1-LZP-thm}). 
 
$\star$ \textit{Subcase 1B. $\deg(u) = 0$, i.e., $u = 1$.}

Note that in this subcase, $P_{u, \fq} = P_{1, \fq}$ since $u = 1$. Recall from Proposition \ref{Proposition-The-1st-basic-results-about-the-Carlitz-annihilator-of-(u, wp)}$(i)$ that $P_{1, \fq}$ divides $\fq - 1$, and since $\gcd(\fq, \fq - 1) = 1$, we deduce that $\gcd(P_{1, \fq}, \fq) = 1$. We should also note from Proposition \ref{Proposition-The-1st-basic-results-about-the-Carlitz-annihilator-of-(u, wp)}$(ii)$ that $\deg(P_{1, \fq}) \ge 1$, and that
\begin{align}
\label{E-inequality-deg(P_{1--fq})--le---deg(fq)---in-Subcase-1B-in-Case-1-LZP-thm}
\deg(P_{1, \fq}) \le \deg(\fq - 1) = \deg(\fq). 
 \end{align}
 
We see from Lemma \ref{Lemma-The-lower-bound-of-Phi(m)} and (\ref{Equation-Eqn-of-m=P_{u-fq}-fq^s--in-Case-1-LZP-the-main-theorem}) that
\begin{align}
\label{E--in--u=1---Phi(m)--ge--s(q-1)^2--deg(P_{1--fq})--deg(fq)----in--Subcase-1B-in-Case1-LZP-thm}
\Phi(m) = \Phi(P_{1, \fq}\fq^s) = \Phi(P_{1, \fq})\Phi(\fq^s) \ge (q - 1)^2\deg(P_{1, \fq})\deg(\fq^s) = s(q - 1)^2\deg(P_{1,\fq})\deg(\fq). 
\end{align}
By Lemma \ref{Lemma-The-degrees-of-cyclotomic-polynomials-when-u-is-a-unit}, we know that
\begin{align}
\label{E-eqn-of-deg(Psi_m(1)--=--dfrac--Phi(m)+delta---q---in-Subcase-1B-in-Case-1-LZP-thm}
\deg(\Psi_m(u)) = \deg(\Psi_m(1)) = \dfrac{\Phi(m) + \delta}{q}
\end{align}
for some integer $\delta \in \{-1, 0, 1\}$. 

Since $u = 1$, we see that $\Psi_m(1) = \Psi_m(u) = \epsilon \fq(m + 1)$. Since $\delta \ge -1$, it thus follows from (\ref{E--in--u=1---Phi(m)--ge--s(q-1)^2--deg(P_{1--fq})--deg(fq)----in--Subcase-1B-in-Case1-LZP-thm}) and (\ref{E-eqn-of-deg(Psi_m(1)--=--dfrac--Phi(m)+delta---q---in-Subcase-1B-in-Case-1-LZP-thm}) that
\begin{align}
\label{E-inequality-of-deg(fq)+deg(m)--ge--s(q-1)^2deg(P_{1--fq})deg(fq)--1--over--q--in--Subcase-1B-in-Case1-LZP-thm}
\deg(\fq) + \deg(m) = \deg(\fq) + \deg(m + 1) = \deg(\Psi_m(1)) \ge \dfrac{s(q - 1)^2\deg(P_{1,\fq})\deg(\fq) - 1}{q}.
\end{align}
By (\ref{E-inequality-deg(P_{1--fq})--le---deg(fq)---in-Subcase-1B-in-Case-1-LZP-thm}), we deduce that
\begin{align*}
\deg(m) = \deg(P_{1, \fq}\fq^s) = \deg(P_{1, \fq}) + s\deg(\fq) \le (s + 1)\deg(\fq), 
\end{align*}
and it thus follows from (\ref{E-inequality-of-deg(fq)+deg(m)--ge--s(q-1)^2deg(P_{1--fq})deg(fq)--1--over--q--in--Subcase-1B-in-Case1-LZP-thm}) that
\begin{align*}
q(s + 2)\deg(\fq) \ge q(\deg(\fq) + \deg(m)) \ge s(q - 1)^2\deg(P_{1,\fq})\deg(\fq) - 1.
\end{align*}
Therefore
\begin{align}
\label{E-the-key-inequality-of-deg(fq)--(s(q-1)^2deg(P-1-fq)--q(s+2))--le--1---in-Subcase-1B-in-Case-1-LZP-thm}
\deg(\fq)(s(q - 1)^2\deg(P_{1,\fq}) - q(s + 2)) \le 1.
\end{align}

We consider the following three subsubcases, according as to whether $\deg(P_{1, \fq}) \ge 3$, $\deg(P_{1, \fq}) = 2$, or $\deg(P_{1, \fq}) = 1$.

$\bullet$ \textit{Subsubcase 1B(i). $\deg(P_{1, \fq}) \ge 3$}

By (\ref{E-inequality-deg(P_{1--fq})--le---deg(fq)---in-Subcase-1B-in-Case-1-LZP-thm}), we see that $\deg(\fq) \ge \deg(P_{1, \fq}) \ge 3$. Since $q \ge 3$ and $s \ge 1$, we deduce that 
\begin{align*}
3sq - (7s + 2) \ge 9s - (7s + 2) = 2s - 2 \ge 0, 
\end{align*}
and thus
\begin{align*}
s(q - 1)^2\deg(P_{1,\fq}) - q(s + 2) \ge 3s(q - 1)^2 - q(s + 2) = q(3sq - (7s + 2)) + 3s \ge 3.
\end{align*}
Therefore $\deg(\fq)(s(q - 1)^2\deg(P_{1,\fq}) - q(s + 2)) \ge 9$, which is a contradiction to (\ref{E-the-key-inequality-of-deg(fq)--(s(q-1)^2deg(P-1-fq)--q(s+2))--le--1---in-Subcase-1B-in-Case-1-LZP-thm}).

$\bullet$ \textit{Subsubcase 1B(ii) $\deg(P_{1, \fq}) = 2$}

We, by appealing to (\ref{E-inequality-deg(P_{1--fq})--le---deg(fq)---in-Subcase-1B-in-Case-1-LZP-thm}), find that 
\begin{align*}
\deg(\fq) \ge \deg(P_{1, \fq}) = 2. 
\end{align*}

If $q \ge 4$, we see that 
\begin{align*}
2sq - (5s + 2) \ge 8s - (5s + 2) = 3s - 2 \ge 1, 
\end{align*}
and thus
\begin{align*}
s(q - 1)^2\deg(P_{1,\fq}) - q(s + 2) = 2s(q - 1)^2 - q(s + 2) = q(2sq - (5s + 2)) + 2s \ge 4 + 2 = 6.
\end{align*}
Hence
\begin{align*}
\deg(\fq)(s(q - 1)^2\deg(P_{1,\fq}) - q(s + 2)) \ge 12,
\end{align*}
which is a contradiction to (\ref{E-the-key-inequality-of-deg(fq)--(s(q-1)^2deg(P-1-fq)--q(s+2))--le--1---in-Subcase-1B-in-Case-1-LZP-thm}). 

Suppose now that $q = 3$. By (\ref{E-the-key-inequality-of-deg(fq)--(s(q-1)^2deg(P-1-fq)--q(s+2))--le--1---in-Subcase-1B-in-Case-1-LZP-thm}), and since $\deg(\fq) \ge \deg(P_{1, \fq}) = 2$, we see that
\begin{align*}
\dfrac{1}{2} \ge \dfrac{1}{\deg(\fq)} \ge s(q - 1)^2\deg(P_{1,\fq}) - q(s + 2) = 8s - 3(s + 2) = 5s - 6,
\end{align*}
and thus
\begin{align*}
s \le \dfrac{13}{10}.
\end{align*}
Since $s$ is a positive integer, we deduce from the last inequality that $s = 1$.

By Proposition \ref{Proposition-The-2nd-result-about-the-degrees-of-cyclotomic-polynomials-when-u-is-a-unit}, we know that
\begin{align*}
\deg(C_{P_{1, \fq}}(1)) = q^{\deg(P_{1, \fq}) - 1} = 3^{2 - 1} = 3,
\end{align*}
and since $C_{P_{1, \fq}}(1) \equiv 0 \pmod{\fq}$ (recall that $P_{1, \fq}$ is the Carlitz annihilator of $(1, \fq)$), and $\deg(\fq) \ge 2$, we deduce that either $\deg(\fq) = 2$ or $\deg(\fq) = 3$. 

In summary, by appealing to (LZP1), we find that the following are true:
\begin{itemize}

\item [(i)] $q = 3$, $u = 1$, and $m = P_{1, \fq}\fq$, where $\fq$ is a monic prime in $\bF_3[T]$ of degree $2$ or $3$ such that the Carlitz annihilator $P_{1, \fq}$ of $(1, \fq)$ is  of degree $2$; 

\item [(ii)] $m + 1$ is a prime in $\bF_3[T]$;

\end{itemize}
This, in particular, implies $m \in \cX_5$, where $\cX_5$ is the set in Lemma \ref{Lemma-The-exceptional-case-EC-V}. Hence $\cX_5 \ne \emptyset$, which is absurd since we know from Lemma \ref{Lemma-The-exceptional-case-EC-V} that $\cX_5 = \emptyset$.

$\bullet$ \textit{Subsubcase 1B(iii). $\deg(P_{1, \fq}) = 1$}

By Proposition \ref{Proposition-The-2nd-result-about-the-degrees-of-cyclotomic-polynomials-when-u-is-a-unit}, we know that
\begin{align*}
\deg(C_{P_{1, \fq}}(1)) = q^{\deg(P_{1, \fq}) - 1} = q^0 = 1,
\end{align*}
and since $C_{P_{1, \fq}}(1) \equiv 0 \pmod{\fq}$ and $\deg(\fq) \ge 1$, we deduce that $\deg(\fq) = 1$. Since $P_{1, \fq}$, $\fq - 1$ are monic polynomials, $P_{1, \fq}$ divides $\fq - 1$ (see Proposition \ref{Proposition-The-1st-basic-results-about-the-Carlitz-annihilator-of-(u, wp)}$(i)$), and $\deg(P_{1, \fq}) = \deg(\fq - 1) = \deg(\fq) = 1$, we deduce that
\begin{align*}
P_{1, \fq} = \fq - 1.
\end{align*}

If $q \ge 5$, we see that 
\begin{align}
\label{E-inequality---sq-(3s+2)---ge---5s--(3s+2)---ge--0--in--Subsubcase-1B(iii)-of-Subcase-1B-in-Case--1-of-LZP-thm}
sq - (3s + 2) \ge 5s - (3s + 2) = 2s - 2 \ge 0.
\end{align}
Note that equality in (\ref{E-inequality---sq-(3s+2)---ge---5s--(3s+2)---ge--0--in--Subsubcase-1B(iii)-of-Subcase-1B-in-Case--1-of-LZP-thm}) occurs if and only if $q = 5$ and $s = 1$.

By (\ref{E-inequality---sq-(3s+2)---ge---5s--(3s+2)---ge--0--in--Subsubcase-1B(iii)-of-Subcase-1B-in-Case--1-of-LZP-thm}), we see that
\begin{align}
\label{E--inequality---deg(fq)-(s(q-1)^2--deg(P--1--fq)--q(s+2))--ge---1----in--Subsubcase-1B(iii)-of-Subcase-1B-in-Case--1-of-LZP-thm}
\deg(\fq)(s(q - 1)^2\deg(P_{1,\fq}) - q(s + 2)) = s(q - 1)^2 - q(s + 2) = q(sq - (3s + 2)) + s \ge 1.
\end{align}
Note that equality in (\ref{E--inequality---deg(fq)-(s(q-1)^2--deg(P--1--fq)--q(s+2))--ge---1----in--Subsubcase-1B(iii)-of-Subcase-1B-in-Case--1-of-LZP-thm}) occurs if and only if $s = 1$ and equality in (\ref{E-inequality---sq-(3s+2)---ge---5s--(3s+2)---ge--0--in--Subsubcase-1B(iii)-of-Subcase-1B-in-Case--1-of-LZP-thm}) occurs.

Combining (\ref{E-the-key-inequality-of-deg(fq)--(s(q-1)^2deg(P-1-fq)--q(s+2))--le--1---in-Subcase-1B-in-Case-1-LZP-thm}) and (\ref{E--inequality---deg(fq)-(s(q-1)^2--deg(P--1--fq)--q(s+2))--ge---1----in--Subsubcase-1B(iii)-of-Subcase-1B-in-Case--1-of-LZP-thm}), we find that 
\begin{align*}
\deg(\fq)(s(q - 1)^2\deg(P_{1,\fq}) - q(s + 2)) = 1,
\end{align*}
which implies that equality in (\ref{E--inequality---deg(fq)-(s(q-1)^2--deg(P--1--fq)--q(s+2))--ge---1----in--Subsubcase-1B(iii)-of-Subcase-1B-in-Case--1-of-LZP-thm}) occurs. By appealing to the remarks following  (\ref{E-inequality---sq-(3s+2)---ge---5s--(3s+2)---ge--0--in--Subsubcase-1B(iii)-of-Subcase-1B-in-Case--1-of-LZP-thm}) and (\ref{E--inequality---deg(fq)-(s(q-1)^2--deg(P--1--fq)--q(s+2))--ge---1----in--Subsubcase-1B(iii)-of-Subcase-1B-in-Case--1-of-LZP-thm}), we find that $q = 5$ and $s = 1$. 

 By appealing to (LZP0) and (LZP1), we find that the following are true:
\begin{itemize}

\item [(i)] $q = 5$, $u = 1$, and $m = P_{1, \fq}\fq^s = (\fq - 1)\fq$, where $\fq$ is a monic prime of degree one in $\bF_5[T]$;

\item [(ii)] $m + 1$ is the only Zsigmondy prime for $(1, m)$; 

\item [(iii)] there are no large Zsigmondy primes for $(1, m)$.

\end{itemize}
This, in particular, implies that $m \in \cX_6$, where $\cX_6$ is the set in Lemma \ref{Lemma-The-exceptional-case-VI}. This is equivalent to saying that we are in the exceptional case (EC-V), which is a contradiction. 

For the rest of \textit{Subsubcase 1B(iii)}, it remains to consider the case when $q = 3$ or $q = 4$.

If $q = 4$, recall that $\deg(\fq) = \deg(P_{1, \fq}) = 1$, and thus we see from (\ref{E-the-key-inequality-of-deg(fq)--(s(q-1)^2deg(P-1-fq)--q(s+2))--le--1---in-Subcase-1B-in-Case-1-LZP-thm}) that
\begin{align*}
1 \ge \deg(\fq)(s(q - 1)^2\deg(P_{1,\fq}) - q(s + 2)) = 5s - 8.
\end{align*}
Hence
\begin{align*}
s \le \dfrac{9}{5},
\end{align*}
and since $s$ is a positive integer, we deduce from the last inequality that $s = 1$. Thus 
\begin{align*}
m = P_{1, \fq}\fq^s = (\fq - 1)\fq,
\end{align*}
where $\fq$ is a monic prime of degree one in $\bF_{2^2}[T]$. This implies that we are in the exceptional case (EC-II), which is a contradiction.

If $q = 3$, we see from (\ref{E-the-key-inequality-of-deg(fq)--(s(q-1)^2deg(P-1-fq)--q(s+2))--le--1---in-Subcase-1B-in-Case-1-LZP-thm}) that
\begin{align*}
1 \ge \deg(\fq)(s(q - 1)^2\deg(P_{1,\fq}) - q(s + 2)) = s - 6,
\end{align*}
and thus 
\begin{align*}
s \le 7.
\end{align*}
Hence 
\begin{align*}
s \in \{1, 2, 3, 4, 5, 6, 7\}. 
\end{align*}
Thus 
\begin{align*}
m = P_{1, \fq}\fq^s = (\fq - 1)\fq^s,
\end{align*}
where $\fq$ is a monic prime of degree one in $\bF_3[T]$ and $s \in \{1, 2, 3, 4, 5, 6, 7\}$. 

If $s = 1$, we deduce that $m = (\fq - 1)\fq$. Since $q = 3$, $u = 1$, and $\fq$ is a monic prime of degree one in $\bF_3[T]$, we see that we are in the exceptional case (EC-I), which is a contradiction. 

If $2 \le s \le 7$, we, by appealing to (LZP0) and (LZP1), find that $m \in \cX_3$, where $\cX_3$ is the set in Lemma \ref{Lemma-The-exceptional-case-EC-III}. This, in turn, is equivalent to saying that we are in the exceptional case (EC-III), which is a contradiction.

$\star$ \textit{Case 2. $\Psi_m(u) = \epsilon (m + 1)$ for some unit $\epsilon \in \bF_q^{\times}$.}

If $\deg(u) \ge 1$, we deduce from Lemma \ref{Lemma-The-degrees-of-cyclotomic-polynomials} and Lemma \ref{Lemma-The-lower-bound-of-Phi(m)} that 
\begin{align*}
\deg(m) = \deg(m + 1) = \deg(\Psi_m(u)) = \Phi(m)\deg(u) \ge (q - 1)\deg(m).
\end{align*}
Thus
\begin{align*}
(q - 2)\deg(m) \le 0,
\end{align*}
which is a contradiction since $\deg(m) \ge 1$ and $q > 2$.

If $\deg(u) = 0$, then $u = 1$, and hence
\begin{align}
\label{E-eqn-of-Psi_m(u)-when-u=1-in-Case-2-of-LZP-thm}
\Psi_m(1) = \epsilon (m + 1).
\end{align}

Since $m$ is of positive degree, there exists a monic prime $\wp$ of positive degree dividing $m$. Then one can write
\begin{align}
\label{E-the-eqn-of-m-in-Case-2-of-LZP-theorem}
m = n\wp^s,
\end{align}
where $s$ is a positive integer, and $n$ is a monic polynomial such that $\gcd(n, \wp) = 1$. 

We consider the following subcases, according as to whether $\deg(n) = 0$ or $\deg(n) \ge 1$.

$\star$ \textit{Subcase 2A. $\deg(n) = 0$.}

In this subcase, since $n$ is monic, we see that $n = 1$, and thus 
\begin{align}
\label{E-eqn-of-m-in-subcase-2A-of-LZP-thm}
m = \wp^s.
\end{align}

We first consider the case when $s = 1$. If $\deg(\wp) = 1$, then $m = \wp$ is a monic prime of degree one in $\bF_q[T]$. This implies that $m \in \cX_7$, where $\cX_7$ is the set in Lemma \ref{L-the-exceptional-case-EC-VI}. Recall that $u = 1$. Hence we are in the exceptional case (EC-VI), which is a contradiction.

Suppose now that $\deg(\wp) \ge 2$. Since $m = \wp$ is a monic prime, we deduce from Lemma \ref{Lemma-The-degrees-of-cyclotomic-polynomials-when-u-is-a-unit} that
\begin{align*}
\deg(\Psi_m(1)) = \deg(\Psi_{\wp}(1)) = \dfrac{\Phi(\wp) + (-1)^2}{q} = \dfrac{\Phi(\wp) + 1}{q},
\end{align*}
and it therefore follows from (\ref{E-eqn-of-Psi_m(u)-when-u=1-in-Case-2-of-LZP-thm}) that
\begin{align}
\label{E-1st-eqn-in-Case-2-of-LZP-thm}
\deg(\wp) = \deg(m) = \deg(\epsilon(m + 1)) = \deg(\Psi_m(1)) = \dfrac{\Phi(\wp) + 1}{q}.
\end{align}

By Corollary \ref{C-lower-bound-for-Phi(prime-of-deg->=2)}, we know that $\Phi(\wp) > q\deg(\wp)$, and thus
\begin{align*}
\dfrac{\Phi(\wp) + 1}{q} > \dfrac{\Phi(\wp)}{q} > \deg(\wp),
\end{align*}
which is a contradiction to (\ref{E-1st-eqn-in-Case-2-of-LZP-thm}). 

We now consider the case when $s \ge 2$. By (\ref{E-eqn-of-Psi_m(u)-when-u=1-in-Case-2-of-LZP-thm}), we know from Lemma \ref{Lemma-The-degrees-of-cyclotomic-polynomials-when-u-is-a-unit} that
\begin{align}
\label{E-1st-eqn-when-s>=2-in-Case-2-of-LZP-thm}
\deg(\wp^s) = \deg(m) =  \deg(\epsilon(m + 1)) = \deg(\Psi_m(1)) = \dfrac{\Phi(\wp^s)}{q},
\end{align}
and it thus follows from Corollary \ref{C-the-weaker-bound-for-Phi(wp^s)} that
\begin{align}
\label{E-2nd-eqn-when-s>=2-in-Case-2-of-LZP-thm}
\deg(\wp^s) =  \dfrac{\Phi(\wp^s)}{q} \ge \dfrac{q\deg(\wp^s)}{q} = \deg(\wp^s).
\end{align}
Therefore equality in (\ref{E-2nd-eqn-when-s>=2-in-Case-2-of-LZP-thm}) occurs, and hence we deduce from Corollary \ref{C-the-weaker-bound-for-Phi(wp^s)} that $q = 3$, $\deg(\wp) = 1$, and $s = 2$. Then $q = 3$, $u = 1$, and $m = \wp^2$, where $\wp$ is a monic prime of degree one in $\bF_3[T]$. By (LZP0), (LZP1), we deduce that $m \in \cX_8$, where $\cX_8$ is the set in Lemma \ref{L-the-exceptional-case-EC-VII}. This implies that we are in the exceptional case (EC-VII), which is a contradiction.

$\star$ \textit{Subcase 2B. $\deg(n)  \ge 1$.}

We first prove that the following is true:
\begin{itemize}

\item [(LZP3)] \textit{$m$ is square-free, that is, $\fq^2$ does not divide $m$ for any monic prime $\fq$.}

\end{itemize}

Assume that (LZP3) does not hold, i.e., there exists a monic prime $\fq$ such that $\fq^2$ divide $m$. Then one can write $m$ in the form
\begin{align}
\label{E-2nd-eqn-of-m-in-terms-of-fq-and-u-in-Subcase-2B-of-LZP-thm}
m = u\fq^r,
\end{align}
where $r \ge 2$, and $u$ is a monic polynomial such that $\gcd(u, \fq) = 1$. Note that since $m = n\wp^s$ (see (\ref{E-the-eqn-of-m-in-Case-2-of-LZP-theorem})), $\gcd(n, \wp) = 1$, $\deg(n) \ge 1$, and $\deg(\wp^s) \ge 1$, we deduce that $\deg(u) > 0$. 

Applying Lemma \ref{Lemma-The-lower-bound-of-Phi(m)} for $u$, and applying Corollary \ref{C-the-weaker-bound-for-Phi(wp^s)} for $\fq^r$, we deduce that
\begin{align}
\label{E-1st-eqn-of-Phi(m)-in-terms-of-Phi(fq^r)-and-Phi(u)-in-Subcase-2B-of-LZP-thm}
\Phi(m) = \Phi(u\fq^r) = \Phi(u)\Phi(\fq^r) \ge (q - 1)\deg(u)q\deg(\fq^r) = q(q - 1)\deg(u)\deg(\fq^r).
\end{align}

We deduce from (\ref{E-eqn-of-Psi_m(u)-when-u=1-in-Case-2-of-LZP-thm}), (\ref{E-2nd-eqn-of-m-in-terms-of-fq-and-u-in-Subcase-2B-of-LZP-thm}), (\ref{E-1st-eqn-of-Phi(m)-in-terms-of-Phi(fq^r)-and-Phi(u)-in-Subcase-2B-of-LZP-thm}), and Lemma \ref{Lemma-The-degrees-of-cyclotomic-polynomials-when-u-is-a-unit} that
\begin{align*}
\deg(u) + \deg(\fq^r) = \deg(u\fq^r) = \deg(m) = \deg(\epsilon(m + 1)) &= \deg(\Psi_m(1)) \\
&= \dfrac{\Phi(m)}{q} \\
&\ge (q - 1)\deg(u)\deg(\fq^r),
\end{align*}
and thus
\begin{align}
\label{E-inequality-to-establish-m-is-not-squarefree-in-Subcase2B-of-LZP-thm}
(q - 2)\deg(u)\deg(\fq^r) + \deg(u)\deg(\fq^r) - (\deg(u) + \deg(\fq^r)) \le 0.
\end{align}

Since $r \ge 2$, we see that $\deg(\fq^r) = r\deg(\fq) \ge 2$, and thus
\begin{align}
\label{E-1st-inequality-of-(q-2)deg(u)deg(fq^r)->=-2-in--m-squarefree-in-Subcase-2B-of-LZP-thm}
(q - 2)\deg(u)\deg(\fq^r) \ge 2(q - 2) \ge 2.
\end{align}
On the other hand, since $\deg(u) \ge 1$ and $\deg(\fq^r) \ge 2$, we know that
\begin{align*}
\deg(u)\deg(\fq^r) - (\deg(u) + \deg(\fq^r)) + 1 = (\deg(u) - 1)(\deg(\fq^r) - 1) \ge 0,
\end{align*}
and thus
\begin{align}
\label{E-2nd-inequality-of-deg(u)deg(fq^r)-(deg(u)-+-deg(fq^r)->=-minus-1-in-Subcase2B-of-LZP-thm}
\deg(u)\deg(\fq^r) - (\deg(u) + \deg(\fq^r)) \ge -1.
\end{align}

From (\ref{E-1st-inequality-of-(q-2)deg(u)deg(fq^r)->=-2-in--m-squarefree-in-Subcase-2B-of-LZP-thm}, (\ref{E-2nd-inequality-of-deg(u)deg(fq^r)-(deg(u)-+-deg(fq^r)->=-minus-1-in-Subcase2B-of-LZP-thm}), we deduce that
\begin{align*}
(q - 2)\deg(u)\deg(\fq^r) + \deg(u)\deg(\fq^r) - (\deg(u) + \deg(\fq^r)) \ge 1,
\end{align*}
which is a contradiction to (\ref{E-inequality-to-establish-m-is-not-squarefree-in-Subcase2B-of-LZP-thm}). Thus (LZP3) is true, i.e., $m$ is square-free. 

In order to get a contradiction in this subcase, we use another representation of $m$. By (\ref{E-the-eqn-of-m-in-Case-2-of-LZP-theorem}), and since $\gcd(n, \wp) = 1$, we deduce that $m$ can be written in the form
\begin{align}
\label{E-eqn-of-m-in-terms-of-m1-and-m2-in-Subcase-2B-of-LZP-thm}
m = m_1m_2,
\end{align}
where $m_1, m_2$ are monic polynomials of positive degrees such that $\gcd(m_1, m_2) = 1$. (For example, one can take $m_1 = n$ and $m_2 = \wp^s$.) Set
\begin{align}
\label{E-eqn-of-alpha-in-Subcase-2B-of-LZP-thm}
\alpha = \deg(m_1) \ge 1,
\end{align}
and
\begin{align}
\label{E-eqn-of-beta-in-Subcase-2B-of-LZP-thm}
\beta = \deg(m_2) \ge 1.
\end{align}
Without loss of generality, one can further assume that
\begin{align}
\label{E-alpha->=-beta-in-Subcase-2B-of-LZP-thm}
\alpha = \deg(m_1) \ge \deg(m_2) = \beta.
\end{align}

Since $\deg(m_1) \ge 1$ and $\deg(m_2) \ge 1$, we deduce from (\ref{E-eqn-of-m-in-terms-of-m1-and-m2-in-Subcase-2B-of-LZP-thm}) and Lemma \ref{Lemma-The-lower-bound-of-Phi(m)} that
\begin{align}
\label{E-inequality-of-Phi(m)-in-terms-of-alpha-and-beta-and-q-in-m=m1m2-in-Subcase-2B-of-LZP-thm}
\Phi(m) = \Phi(m_1m_2) = \Phi(m_1)\Phi(m_2) \ge (q - 1)^2\deg(m_1)\deg(m_2) = (q - 1)^2\alpha\beta.
\end{align}

Since $m$ is square-free, we deduce from  (\ref{E-eqn-of-Psi_m(u)-when-u=1-in-Case-2-of-LZP-thm}) and Lemma \ref{Lemma-The-degrees-of-cyclotomic-polynomials-when-u-is-a-unit} that there exists an integer $\delta \in \{-1, 1\}$ such that
\begin{align}
\label{E-alpha-+-beta-equals-the-quotient-of-Phi(m)+delta-over-q-in-m=m1m2-in-Subcase-2B-of-LZP-thm}
\alpha + \beta = \deg(m_1) + \deg(m_2) = \deg(m) = \deg(\epsilon(m + 1)) = \deg(\Psi_m(1)) = \dfrac{\Phi(m) + \delta}{q}.
\end{align}
Since $\delta \ge -1$, we deduce from (\ref{E-inequality-of-Phi(m)-in-terms-of-alpha-and-beta-and-q-in-m=m1m2-in-Subcase-2B-of-LZP-thm}) and (\ref{E-alpha-+-beta-equals-the-quotient-of-Phi(m)+delta-over-q-in-m=m1m2-in-Subcase-2B-of-LZP-thm}) that
\begin{align}
\label{E-alpha-+-beta-ge-(q-1)^2-alpha-beta-minus-1-over-q-in-m=m1m2-in-Subcase-2B-of-LZP-thm}
\alpha + \beta \ge \dfrac{(q - 1)^2\alpha\beta - 1}{q}.
\end{align}
Inequality (\ref{E-alpha-+-beta-ge-(q-1)^2-alpha-beta-minus-1-over-q-in-m=m1m2-in-Subcase-2B-of-LZP-thm}) is equivalent to the inequality
\begin{align}
\label{E-the-main-inequality-involving-q-alpha-beta-le-1--in-m=m1m2--in-Subcase-2B-of-LZP-thm}
q((\alpha\beta)q - (2\alpha\beta + \alpha + \beta)) + \alpha\beta = (\alpha\beta)q^2 - (2\alpha\beta + \alpha + \beta)q + \alpha\beta \le 1.
\end{align}

Note that equality in (\ref{E-alpha-+-beta-ge-(q-1)^2-alpha-beta-minus-1-over-q-in-m=m1m2-in-Subcase-2B-of-LZP-thm}) occurs if and only if $\delta = -1$ and equality in (\ref{E-inequality-of-Phi(m)-in-terms-of-alpha-and-beta-and-q-in-m=m1m2-in-Subcase-2B-of-LZP-thm}) occurs. Since inequality (\ref{E-alpha-+-beta-ge-(q-1)^2-alpha-beta-minus-1-over-q-in-m=m1m2-in-Subcase-2B-of-LZP-thm}) is equivalent to inequality (\ref{E-the-main-inequality-involving-q-alpha-beta-le-1--in-m=m1m2--in-Subcase-2B-of-LZP-thm}), we should also note that the following is true:
\begin{itemize}

\item [(ELZP)] \textit{Equality in (\ref{E-the-main-inequality-involving-q-alpha-beta-le-1--in-m=m1m2--in-Subcase-2B-of-LZP-thm}) occurs if and only if $\delta = -1$ and equality in (\ref{E-inequality-of-Phi(m)-in-terms-of-alpha-and-beta-and-q-in-m=m1m2-in-Subcase-2B-of-LZP-thm}) occurs}.

\end{itemize}

Since $\alpha \ge 1$ and $\beta \ge 1$, we deduce that 
\begin{align*}
\alpha\beta - (\alpha + \beta) + 1 = (\alpha - 1)(\beta - 1) \ge 0, 
\end{align*}
and thus
\begin{align}
\label{E-inequality-of-alpha-beta--alpha-beta-ge-alpha-plus-beta-minus-1--in-m=m1m2--in-Subcase-2B-of-LZP-thm}
\alpha\beta \ge \alpha + \beta - 1.
\end{align}
Note that equality in (\ref{E-inequality-of-alpha-beta--alpha-beta-ge-alpha-plus-beta-minus-1--in-m=m1m2--in-Subcase-2B-of-LZP-thm}) occurs if and only if $\alpha = 1$ or $\beta = 1$. 

Since $q \ge 3$, we deduce from (\ref{E-inequality-of-alpha-beta--alpha-beta-ge-alpha-plus-beta-minus-1--in-m=m1m2--in-Subcase-2B-of-LZP-thm}) that
\begin{align}
\label{E-inequality-of-alpha-beta-and-q-----(alpha-beta)q-minus-(2alpha-beta-plus-alpha-plus-beta)-ge-minus-1--in-m=m1m2--in-Subcase-2B-of-LZP-thm}
(\alpha\beta)q - (2\alpha\beta + \alpha + \beta) \ge 3\alpha\beta - (2\alpha\beta + \alpha + \beta) = \alpha\beta - (\alpha + \beta) \ge - 1.
\end{align}
Since $(\alpha\beta)q - (2\alpha\beta + \alpha + \beta)$ is an integer, we deduce from the above inequality that either 
\begin{align}
\label{E-the-1st-main-eqn-of-alpha-beta-and-q-ge-0--in-m=m1m2--in-Subcase-2B-of-LZP-thm}
(\alpha\beta)q - (2\alpha\beta + \alpha + \beta) \ge 0
\end{align}
or 
\begin{align}
\label{E-the-2nd-main-eqn-of-alpha-and-beta-and-q-eq-minus-1--in-m=m1m2--in-Subcase-2B-of-LZP-thm}
(\alpha\beta)q - (2\alpha\beta + \alpha + \beta) = -1.
\end{align}
Note that (\ref{E-the-2nd-main-eqn-of-alpha-and-beta-and-q-eq-minus-1--in-m=m1m2--in-Subcase-2B-of-LZP-thm}) holds if and only if equalities in (\ref{E-inequality-of-alpha-beta-and-q-----(alpha-beta)q-minus-(2alpha-beta-plus-alpha-plus-beta)-ge-minus-1--in-m=m1m2--in-Subcase-2B-of-LZP-thm}) occur at the same time. This of course implies that (\ref{E-the-2nd-main-eqn-of-alpha-and-beta-and-q-eq-minus-1--in-m=m1m2--in-Subcase-2B-of-LZP-thm}) holds if and only if $q = 3$, and $\alpha = 1$ or $\beta = 1$. 

We consider the following subsubcases:

$\bullet$ \textit{Subsubcase 2B(i). (\ref{E-the-1st-main-eqn-of-alpha-beta-and-q-ge-0--in-m=m1m2--in-Subcase-2B-of-LZP-thm}) holds, i.e., $(\alpha\beta)q - (2\alpha\beta + \alpha + \beta) \ge 0$.}

In this subsubcase, since $\alpha \ge 1$ and $\beta \ge 1$, we see that
\begin{align}
\label{E-1st-inequality-of-alpha-beta-and-q-in-m=m1m2--in-Subsubcase-2B(i)-of-LZP-thm}
q((\alpha\beta)q - (2\alpha\beta + \alpha + \beta)) + \alpha\beta \ge 1.
\end{align}
Note that equality in (\ref{E-1st-inequality-of-alpha-beta-and-q-in-m=m1m2--in-Subsubcase-2B(i)-of-LZP-thm}) occurs if and only 
\begin{align}
\label{E-2nd-eqn-of-alpha-beta-and-q-in-m=m1m2-in-Subsubcase-2B(i)-of-LZP-thm}
(\alpha\beta)q - (2\alpha\beta + \alpha + \beta) = 0,
\end{align}
and $\alpha = \beta = 1$.

From (\ref{E-the-main-inequality-involving-q-alpha-beta-le-1--in-m=m1m2--in-Subcase-2B-of-LZP-thm}) and (\ref{E-1st-inequality-of-alpha-beta-and-q-in-m=m1m2--in-Subsubcase-2B(i)-of-LZP-thm}), we deduce that equality in (\ref{E-1st-inequality-of-alpha-beta-and-q-in-m=m1m2--in-Subsubcase-2B(i)-of-LZP-thm}) occurs, and it thus follows from the above remark that $\alpha = \beta = 1$, and 
\begin{align}
\label{E-3rd-eqn-of-alpha-beta-and-q-implies-that-q=4--in-m=m1m2-in-Subsubcase-2B(i)-of-LZP-thm}
(\alpha\beta)q - (2\alpha\beta + \alpha + \beta) = 0.
\end{align}
Since $\alpha = \beta = 1$, equation (\ref{E-3rd-eqn-of-alpha-beta-and-q-implies-that-q=4--in-m=m1m2-in-Subsubcase-2B(i)-of-LZP-thm}) implies that $q = 4$. 

In summary, we have showed in this subsubcase that $q = 4$, $u = 1$, and $m = m_1m_2$, where $m_1, m_2$ are monic polynomials in $\bF_4[T]$ such that $\deg(m_1) = \deg(m_2) = 1$ and $\gcd(m_1, m_2) = 1$. It then follows from (LZP0) and (LZP1) that $m \in \cX_9$, where $\cX_9$ is the set in Lemma \ref{L-the-exceptional-case-EC-VIII}. This implies that we are in the exceptional case (EC-VIII), which is a contradiction.

$\bullet$ \textit{Subsubcase 2B(ii). (\ref{E-the-2nd-main-eqn-of-alpha-and-beta-and-q-eq-minus-1--in-m=m1m2--in-Subcase-2B-of-LZP-thm}) holds, i.e., $(\alpha\beta)q - (2\alpha\beta + \alpha + \beta) = -1$.}

The remark following (\ref{E-the-2nd-main-eqn-of-alpha-and-beta-and-q-eq-minus-1--in-m=m1m2--in-Subcase-2B-of-LZP-thm}) tells us that in this subsubcase, $q = 3$, and $\alpha = 1$ or $\beta = 1$. If $\alpha = 1$, we see from (\ref{E-alpha->=-beta-in-Subcase-2B-of-LZP-thm}) that $\beta = 1$. Thus, in any event, $q = 3$ and $\beta = 1$. 

Since $q = 3$ and $\beta = 1$, we deduce from (\ref{E-the-main-inequality-involving-q-alpha-beta-le-1--in-m=m1m2--in-Subcase-2B-of-LZP-thm}) and  (\ref{E-the-2nd-main-eqn-of-alpha-and-beta-and-q-eq-minus-1--in-m=m1m2--in-Subcase-2B-of-LZP-thm}) that $\alpha \le 4$, and therefore $\alpha \in \{1, 2, 3, 4\}$. 

We contend that $\alpha = 1$ or $\alpha = 2$. Indeed, if $\alpha = 4$, then one sees that equality in (\ref{E-the-main-inequality-involving-q-alpha-beta-le-1--in-m=m1m2--in-Subcase-2B-of-LZP-thm}) occurs, and it thus follows from (ELZP) that equality in (\ref{E-inequality-of-Phi(m)-in-terms-of-alpha-and-beta-and-q-in-m=m1m2-in-Subcase-2B-of-LZP-thm}) occurs. This implies that 
\begin{align}
\label{E-1st-eqn-of-Phi(m)-equals-16--in-Subsubcase-2B(ii)-of-LZP-thm}
\Phi(m) = \Phi(m_1)\Phi(m_2) = (q - 1)^2\alpha\beta = 16.
\end{align}
Since $\deg(m_2) = \beta = 1$, we see that $m_2$ is a monic prime of degree one in $\bF_3[T]$. Hence $\Phi(m_2) = 3^{\deg(m_2)} - 1 = 3 - 1 = 2$, and thus
\begin{align}
\label{E-2nd-eqn-of-Phi(m1)-equals-8-in-Subsubcase-2B(ii)-of-LZP-thm}
\Phi(m_1) = 8.
\end{align}

We should note that $m_1$ is square-free since $m$ is square-free. Since $\alpha = \deg(m_1) = 4$, either all monic prime factors of $m_1$ are of degree one or there exists a monic prime, say $P$, of degree at least 2 such that $P$ divides $m_1$. If the former holds, then $m_1$ is of the form
\begin{align*}
m_1 = P_1P_2P_3P_4,
\end{align*}
where the $P_i$ are distinct monic primes of degree one in $\bF_3[T]$. Then
\begin{align*}
\Phi(m_1) = \Phi(P_1)\Phi(P_2)\Phi(P_3)\Phi(P_4) = (3^{\deg(P_1)} - 1)(3^{\deg(P_2)} - 1)(3^{\deg(P_3)} - 1)(3^{\deg(P_4)} - 1) = 16,
\end{align*}
which is a contradiction to (\ref{E-2nd-eqn-of-Phi(m1)-equals-8-in-Subsubcase-2B(ii)-of-LZP-thm}). 

Suppose now that there exists a monic prime $P$ of degree at least 2 such that $P$ divides $m_1$. Write
\begin{align*}
m_1 = Pm_{*},
\end{align*}
where $m_{*}$ is a monic polynomial such that $\gcd(m_{*}, P) = 1$. Since $\deg(m_1) = \alpha = 4$, we see that $2 \le \deg(P) \le 4$.

If $\deg(P) = 4$, then $m_1 = P$, and thus 
\begin{align*}
\Phi(m_1) = \Phi(P) = 3^{\deg(P)} - 1 = 3^4 - 1 = 80,
\end{align*}
which is a contradiction to (\ref{E-2nd-eqn-of-Phi(m1)-equals-8-in-Subsubcase-2B(ii)-of-LZP-thm}).

If $2 \deg(P) \le 3$, then $\deg(m_{*}) \ge 1$. We see from Lemma \ref{Lemma-The-lower-bound-of-Phi(m)} that 
\begin{align*}
\Phi(m_{*}) \ge (q - 1)\deg(m_{*}) \ge 2.
\end{align*}
On the other hand, since $2 \deg(P) \le 3$, we deduce that 
\begin{align*}
\Phi(P) = 3^{\deg(P)} - 1 \ge 3^2 - 1 = 8.
\end{align*}
Therefore
\begin{align*}
\Phi(m) = \Phi(P)\Phi(m_{*}) \ge 16,
\end{align*}
which is a contradiction to (\ref{E-2nd-eqn-of-Phi(m1)-equals-8-in-Subsubcase-2B(ii)-of-LZP-thm}). Thus, by what we have showed above, we deduce that $\alpha \in \{1, 2, 3\}$. 

Suppose now that $\alpha = 3$. Recall that $q = 3$ and $\beta = 1$. We now deduce from (\ref{E-alpha-+-beta-equals-the-quotient-of-Phi(m)+delta-over-q-in-m=m1m2-in-Subcase-2B-of-LZP-thm}) that
\begin{align}
\label{E-eqn-of-Phi(m)-equals-a-prime-when-alpha=3-in-Subsubcase-2B(ii)-of-LZP-thm}
4 = \alpha + \beta = \dfrac{\Phi(m) + \delta}{q} = \dfrac{\Phi(m) + \delta}{3}.
\end{align}
Since $\delta = \pm 1$, we find from (\ref{E-eqn-of-Phi(m)-equals-a-prime-when-alpha=3-in-Subsubcase-2B(ii)-of-LZP-thm}) that either $\Phi(m) = 11$ or $\Phi(m) = 13$. On the other hand, since $\beta = \deg(m_2) = 1$, we deduce that $m_1$ is a monic prime of degree one, and it thus follows from (\ref{E-eqn-of-m-in-terms-of-m1-and-m2-in-Subcase-2B-of-LZP-thm}) that
\begin{align*}
\Phi(m) = \Phi(m_1)\Phi(m_2) = \Phi(m_1)(3^{\deg(m_2)} - 1) = 2\Phi(m_1) \equiv 0 \pmod{2},
\end{align*}
which is a contradiction to the fact that either $\Phi(m) = 11$ or $\Phi(m) = 13$. This contradiction implies that $\alpha = 1$ or $\alpha = 2$. 

In summary, we have showed in \textit{Subsubcase 2B(ii)} that $q = 3$, $u = 1$, and $m = m_1m_2$, where $m_1, m_2$ are monic polynomials in $\bF_3[T]$ such that $\gcd(m_1, m_2) = 1$, $\alpha = \deg(m_1) \in \{1, 2\}$, and $\beta = \deg(m_2) = 1$. We therefore, by appealing to (LZP0), (LZP1), and (LZP3), find that $m \in \cX_{10}$, where $\cX_{10}$ is the set in Lemma \ref{L-the-exceptional-case-EC-IX}. This implies that we are in the exceptional case (EC-IX), which is a contradiction. 

In any event, by all of what we have showed above, we deduce that there exists a large Zsigmondy prime for $(u, m)$, and therefore Theorem \ref{Theorem-The-main-theorem-about-large-Zsigmondy-primes} follows immediately.

\end{proof}

\section*{Acknowledgements}

I would like to thank the referee very much for a very careful reading, pointing out some relevant and useful references, and making many thoughtful comments on an earlier version of my paper. All the calculations in this paper were performed using the computational algebra software MAGMA (see \cite{magma}).

\end{document}